\documentclass[twocolumn]{autart}

\usepackage[T1]{fontenc} 
\usepackage{cite}

\usepackage{braket,amsfonts}
\usepackage{graphicx}
\usepackage{amsmath}
\usepackage{amssymb}
\usepackage{amsfonts}
\usepackage{hyperref}
\usepackage[caption=false]{subfig}
\usepackage{pgfplots}
\usepackage{algorithm}
\usepackage[noend]{algpseudocode}
\usepackage{graphics}
\usepackage{color}
\usepackage{epsfig}
\usepackage{float}
\usepackage{wrapfig}
\usepackage{bbold}
\usepackage{appendix}
\usepackage{hyperref}
\usepackage{textcomp}
\usepackage{graphicx,color}
\usepackage{mathrsfs}
\usepackage{url}
\usepackage{color}
\usepackage{dsfont}
\usepackage{bbm}
\usepackage{verbatim}
\usepackage{tikz}

\newtheorem{theorem}{Theorem}
\newtheorem{lemma}{Lemma}

\newtheorem{remark}{Remark}
\newtheorem{assumption}{Assumption}

\newcommand{\real}{{\mathbb{R}}}

\newcommand{\realnonnegative}{{\mathbb{R}}_{\ge 0}}

\newcommand{\naturals}{\mathbb{N}}

\newcommand{\dvol}{\operatorname{dvol}}
\newcommand{\dd}{\operatorname{d}}
\newcommand{\vol}{\operatorname{vol}}

\newcommand{\Lip}{\operatorname{Lip}}

\newcommand\oprocendsymbol{\hbox{$\bullet$}}
\newcommand\oprocend{\relax\ifmmode\else\unskip\hfill\fi\oprocendsymbol} 

\makeatletter
\renewcommand*{\@opargbegintheorem}[3]{\trivlist
      \item[\hskip \labelsep{\bfseries #1\ #2}] \textbf{(#3)}\ \itshape}
\makeatother    

\allowdisplaybreaks

\begin{document}

\begin{frontmatter}

\title{Distributed Online Optimization for \\ Multi-Agent Optimal Transport
		}

\thanks[footnoteinfo]{This material is based upon work supported by grant AFOSR FA9550-18-1-0158 and ONR N00014-19-12471.}

\author[krishnan]{Vishaal Krishnan} \ead{vkrishnan@seas.harvard.edu},      
\author[martinez]{Sonia Mart\'inez} \ead{soniamd@ucsd.edu}               

\address[krishnan]{School of Engineering and Applied Sciences, Harvard University, Cambridge MA 02138 USA.}                               
\address[martinez]{Department of Mechanical and Aerospace Engineering, University of California, San Diego, La Jolla CA 92093 USA.}             


\begin{abstract}
  We propose a scalable, distributed algorithm for the
  optimal transport of large-scale multi-agent systems.  
  We formulate the problem as one of steering 
  the collective towards a target probability measure while minimizing
  the total cost of transport, with the additional constraint of
  distributed implementation. Using optimal transport theory, 
  we realize the solution as an iterative transport 
  based on a stochastic proximal descent scheme.
  At each stage of the transport, the agents 
  implement an online, distributed primal-dual algorithm 
  to obtain local estimates of the Kantorovich potential
  for optimal transport from the current distribution of the collective 
  to the target distribution.
  Using these estimates as their local objective functions,
  the agents then implement the transport by stochastic proximal descent. 
  This two-step process is carried out recursively by the
  agents to converge asymptotically to the target distribution.  
  We rigorously establish the underlying theoretical framework and convergence of the algorithm and test its behavior in numerical experiments.
\end{abstract}

\end{frontmatter}

\section{Introduction}
We consider the problem of designing
distributed feedback control laws for optimally transporting 
a group of agents from an initial configuration, given by
a probability measure~$\mu_0$, to a target configuration
specified by a probability measure~$\mu^*$ over a domain~$\Omega$.
Such a problem of transport of multi-agent collectives arises 
naturally in various settings, from the modeling of
cell populations in biology, to engineering applications of
coverage control and deployment in robotics and mobile sensing 
networks~\cite{FB-JC-SM:09, MM-ME:10, JC-SM-TK-FB:02j}.
As these scenarios involve physical transport of resources, there is an associated cost of transport owing to energy considerations.
Optimal transport theory~\cite{CV:08, YC-TG-MP:21}, which deals with the problem of rearranging probability measures while minimizing
a cost of transport, presents an appropriate theoretical framework
for the problem. Another consideration in the multi-agent setting
is the scalability of implementation when the size of the collective 
increases, which underlines the need for distributed algorithms~\cite{FB-JC-SM:09}.
This need is further exacerbated by limitations on sensing and
communication typically present in these systems, whereby individual agents may
only be able to sense/communicate with their spatial neighbors.
However, existing formulations of the optimal transport problem 
result in solutions wherein the optimal transport routes for the individual agents depend on the distribution of the entire multi-agent system, which is a bottleneck for distributed implementation.
The primary challenge in algorithm design, therefore, is in obtaining scalable distributed algorithms for large-scale implementation of optimal transport, which constitutes the aim of this paper.

\textbf{Related work.}
Transport problems in robotics and mobile sensing network applications arise in the
form of coverage control and deployment objectives, where the 
underlying goal is to steer a group of robots
towards a target coverage profile over a spatial region.
Among the approaches to the coverage control and deployment problem
for large-scale multi-agent systems are 
transport by synthesis of Markov 
transition matrices~\cite{SB-SJC-FYH:13, ND-UE-BA:15,
SB-SJC-FH:17}, the use of continuum models~\cite{VK-SM:18-sicon, UE-BA:17, TZ-ZL-HL:20} for transport, 
and coverage control by parameter tuning and/or boundary control of the
reaction-advection-diffusion PDE~\cite{KE-HK-SB:18,FZ-AB-KE-SB:18}.
We note, however, that despite the potential for the
application of optimal transport ideas to the multi-agent setting,
as seen from the works~\cite{SF-GF-PZ-TAW:16, PZ-CL-SF:21, SB-SJC-FYH:14, CS-SB-MQ-FB:21},
a truly distributed formulation of optimal transport has remained open.
While our preliminary work~\cite{VK-SM:18-cdc} was an attempt in this direction,
we develop a rigorous theoretical framework for the design of distributed 
optimal transport algorithms in this paper.
These works, however, present significant limitations
either because they require centralized offline planning~\cite{SF-GF-PZ-TAW:16}, 
or because of a need for costly computation and information
exchange between agents~\cite{SB-SJC-FYH:14}.
This serves as a strong motivation for the 
development of a distributed iterative algorithm for optimal transport
in this paper.

The applications of optimal transport in image processing and 
various engineering domains has motivated a search for efficient
computational methods for the optimal transport 
problem~\cite{MC:13, GP-MC:17, SL-SM-YC-HZ-HZ:21}.
Optimal transport from continuous to discrete probability distributions
has been studied under the name of semi-discrete optimal transport,
with connections to the problem of optimal quantization of
probability measures, in~\cite{DB-BS-BW:18}.
While computational approaches to optimal transport
often work with the static, Monge or Kantorovich formulations
of the problem, investigations involving dynamical formulations was initiated 
by~\cite{JDB-YB:00}, where the authors recast the~$L^2$ 
Monge-Kantorovich mass transfer problem in a fluid mechanics 
framework. This largely owes to notion of displacement interpolation 
originally introduced in~\cite{RM:95}. 
\cite{NP-GP-EO:14} and~\cite{JDB-GC-MC-LN-GP:15} are other 
works in this vein. The problem of optimal transport was also 
explored from a stochastic control perspective in~\cite{TM-MT:08, YC-TG-MP:16, MHDB-EM-RS:21}.
However, there has remained a gap in this literature  
with regard to distributed computation of optimal transport,
which arises as a rather stringent constraint in multi-agent transport scenarios.


\textbf{Contributions.}
In this work, we propose and investigate  
large-scale optimal transport of multi-agent
collectives based on a scalable, distributed online optimization.
Working with a reduction of the Kantorovich duality for metric costs
conformal to the Euclidean metric, we note that the Kantorovich potential
is almost everywhere differentiable and obtain a bound on the norm of its
gradient. We then obtain an stochastic process for optimal iterative transport of probability 
measures based on Kantorovich duality, showing it to be equivalent 
to optimal transport along geodesics,
and establish convergence of the sequence of
probability measures generated by the process 
to the target probability measure with respect to 
the topology of weak convergence.
We propose a distributed primal-dual algorithm 
to be implemented online by the agents to obtain
local estimates of the Kantorovich potential, which are then
used as local objectives in a proximal algorithm for transport.
The paper contributes not only to the literature on computational
methods for the optimal transport problem,
but also presents a novel scalable, distributed approach to
multi-agent optimal transport addressing a longstanding 
concern in the research on multi-agent systems.

\textbf{Paper outline.} 
In Section~\ref{sec:prelims}, we introduce the notation and mathematical
preliminaries underlying the results presented in the rest of the paper.
Our goal in this paper is to design an iterative transport process 
relying entirely on distributed computation
for the optimal transport of a multi-agent system 
towards a target probability measure.
To this end, we first design in Section~\ref{sec:iterative_OT} 
an iterative process for optimally transporting a 
probability measure onto a target measure with respect to an underlying  transport cost on the spatial domain. 
We then obtain in Section~\ref{sec:DOT} a distributed multi-agent optimal transport algorithm via a discretization of the transport process introduced
in Section~\ref{sec:iterative_OT}, which is then followed by
an investigation of the behavior of the algorithm in 
numerical simulations. We then conclude with a brief summary of
our results in Section~\ref{sec:conclusion}.

\section{Mathematical preliminaries} \label{sec:prelims}
\subsection{Notation}
We first briefly introduce the notation adopted in the rest of the paper.
We use~$| \cdot |$ to denote the Euclidean norm in~$\real^d$, for
any~$d \in \mathbb{N}$ (when $d=1$, this denotes the absolute value).
We use~$\| \cdot \|$ for function space norms.
The gradient operator in $\real^d$ is represented as $\nabla$.
  For any $\Omega \subseteq \real^d$,~$\partial \Omega \subseteq \real^d$ denotes its boundary,
  $\bar{\Omega} = \Omega \cup \partial \Omega$ its closure, and
  $\mathring{\Omega} = \Omega \setminus \partial \Omega$ its interior
  with respect to the standard Euclidean topology. 
  We denote by $(\Omega, \mu)$ the set~$\Omega$ with an underlying measure~$\mu$.
  Given any~$x \in \Omega \subset \real^d$, the set~$B_r(x)$ is the closed $d$-ball of
  radius~$r>0$, centered at~$x$. Furthermore, we denote by $B^c_{r}(x)$ the closed $d$-ball of radius $r$, centered at~$x$, with respect to the metric~$c$.
  We denote by $\mathcal{P}(\Omega)$ the space of probability measures
  over $\Omega$. For a measurable mapping~$\mathcal{T} : \Omega
  \rightarrow \Theta$, where $\Omega$ and $\Theta$ are measurable, we
  denote by $\mathcal{T}_{\#} \mu \in \mathcal{P}(\Theta)$ the
  pushforward measure of~$\mu \in \mathcal{P}(\Omega)$ and we have
  $\mathcal{T}_{\#} \mu (B) = \mu (\mathcal{T}^{-1}(B))$,  for all
  measurable $B \subseteq \Theta$.
  We use~$\left \langle f, g \right \rangle$ to represent the inner
  product of functions~$f,g : \Omega \rightarrow \real$ w.r.t.~the
  Lebesgue measure~$\vol$, given by~$\left \langle f,g \right \rangle =
  \int_{\Omega} fg \dvol$. The set $\Lip(\Omega)$ is the space of
  Lipschitz continuous functions on $\Omega$.
We denote by~$L^p(\Omega, \mu)$ the space of~$p$-integrable (measurable) 
functions on~$\Omega$, where the integration is carried out with the underlying
measure~$\mu$ (the Lebesgue measure is implied when~$\mu$ is not specified),
and by~$W^{1,p}(\Omega, \mu)$ the space of $p$-integrable (measurable) 
functions with $p$-integrable (measurable) derivatives.

\subsection{Monge and Kantorovich formulations of optimal transport}
Let~$\Omega \subseteq \real^d$ be a compact, convex domain, 
and let~$\mu, \nu \in \mathcal{P}(\Omega)$ be absolutely continuous
probability measures on~$\Omega$.  Let~$c : \Omega \times \Omega
\rightarrow \realnonnegative$ be a continuous function 
such that for~$x,y \in \Omega$, $c(x,y)$ is the unit cost of transport from~$x$ to~$y$. 
In the Monge (deterministic) formulation,
the optimal cost of transporting the probability measure~$\mu$ 
onto~$\nu$ is defined as the infimum of the transport cost
over the set of maps for which~$\nu$ is 
obtained as the pushforward measure of~$\mu$,
as given below:
\begin{align}
  C_M(\mu, \nu) = \inf_{ \substack{T:\Omega \rightarrow \Omega \\
      T_{\#}\mu = \nu }} ~\int_{\Omega} c(x,T(x))
  d\mu(x).
  	\label{eq:OT_Monge}
\end{align}
The Kantorovich (probabilistic) formulation relaxes the 
Monge formulation~\eqref{eq:OT_Monge}
by defining the optimal cost of transporting the probability measure~$\mu$ 
onto~$\nu$ as the infimum of the transport cost over the set of joint probability 
measures~$\Pi(\mu, \nu) \subset \mathcal{P}(\Omega \times \Omega)$,
for which~$\mu$ and~$\nu$ are the respective marginals over~$\Omega$,
as given below:
\begin{align}
  C_K(\mu, \nu) = \inf_{\pi \in \Pi(\mu, \nu)} \int_{\Omega \times
    \Omega} c(x,y)~d\pi(x,y).
	\label{eq:OT_Kantorovich}
\end{align}
Since~$\Omega$ is a compact subset of~$\real^d$ and
$c$ is continuous over $\Omega \times \Omega$,
it follows that the Kantorovich problem admits a solution~\cite{FS:15}. 
Furthermore, the Kantorovich formulation can be shown~\cite{FS:15} to be a 
relaxation of the Monge formulation, i.e. $C_M(\mu, \nu) = C_K(\mu, \nu)$,  
and we hereafter denote by~$C(\mu, \nu) = C_M(\mu, \nu) = C_K(\mu, \nu)$. 
Furthermore, the Kantorovich formulation~\eqref{eq:OT_Kantorovich} admits the
following dual formulation\footnote{Strong duality holds for the 
Kantorovich formulation (c.f.~Theorem~5.10 in~\cite{CV:08}).}:
\begin{align}
	\begin{aligned}
 		 C(\mu, \nu) = 
 		 \sup_{\substack{\phi, \psi \in L^1(\Omega)}}
  		~&\int_{\Omega} \phi(x) d\mu(x) + \int_{\Omega} \psi(y) d\nu(y) \\
  		&\text{s.t~~} \phi(x) + \psi(y) \leq c(x,y).
	 \end{aligned}
  \label{eq:Kantorovich_duality}
\end{align}
The maximizers of the above dual formulation are pairs of 
functions~$(\phi, \psi)$ called Kantorovich potentials. They 
occur at the boundary of the inequality constraint, thereby satisfying:
\begin{align}
	\begin{aligned}
          \phi(x) &= \inf_{y \in \Omega} \left( c(x,y) - \psi(y) \right), \\
          \psi(y) &= \inf_{z \in \Omega} \left( c(z,y) - \phi(z) \right).
	\end{aligned}
	\label{eq:conjugate_pair}
\end{align}
We refer to~$(\phi, \psi)$ defined above as a~$c$-conjugate pair, and
write~$\psi = \phi^c$ to denote that $\psi$ is the conjugate of $\phi$. 

\section{Optimal iterative transport of measures}
\label{sec:iterative_OT}

We first briefly describe the setting for the multi-agent optimal
transport problem addressed in this paper.
We consider a compact and convex domain $\Omega \subset \real^d$
across which $N$~agents are initially independently and
identically distributed according to an absolutely
continuous probability measure~$\mu_0 \in \mathcal{P}(\Omega)$,
i.e., the initial agent positions 
$\left \lbrace x_1(0), \ldots, x_N(0) \right \rbrace$
are independently generated as $x_i(0) \sim \mu_0$ for any $i \in \lbrace 1, \ldots, N \rbrace$. 
Our goal in this paper is to design an iterative transport process 
to steer the agents towards an absolutely continuous target probability measure $\mu^* \in \mathcal{P}(\Omega)$, while minimizing the net cost of transport
measured with respect to an underlying metric~$c$ on~$\Omega$.
In particular, we restrict the transport cost~$c$
to the class of distance functions induced by a Riemannian metric conformal to the Euclidean one,
as stated in the following assumption:
\begin{assumption}[\bf \emph{Conformal distance}] \label{ass:conformal}
Let~$c: \Omega \times \Omega \rightarrow \realnonnegative$ 
be a distance function on~$\Omega$ conformal to the Euclidean distance
(with a strictly positive conformal factor~$\xi \in C^1(\Omega)$),
i.e., for any~$x,y \in \Omega$, $c(x,y)$ is given by:
\begin{align} \label{eq:conformal_metric}
	c(x,y) = &\inf_{\substack{\gamma \in C^1 \left([0,1];  \Omega \right); \\ \gamma(0) = x, \gamma(1) = y}}  
									 ~\int_0^1 \xi(\gamma(t)) \; | \; \dot{\gamma}(t) | \dd t.
\end{align}
\end{assumption}
The restriction to a conformal metric affords us the flexibility to penalize transport through certain regions of the domain relative to
other regions via the conformal factor, therefore allowing for the possibility of modeling physical obstacles and uneven terrain in real environments with relative ease.
The following lemma\footnote{We refer the reader to the extended version~\cite{VK-SM:18-arxiv} for the detailed proofs of the results in the paper.} establishes a reduction of the 
Kantorovich duality~\eqref{eq:Kantorovich_duality}
for conformal metric costs:
\begin{lemma}[\bf \emph{Local Lipschitz bound}]
\label{lemma:cconvex_phic}
Let $c: \Omega \times \Omega \rightarrow \realnonnegative$ 
be a metric on~$\Omega$ conformal to the Euclidean distance, 
with a strictly positive conformal factor~$\xi \in C^1(\Omega)$.
(a) The conjugate of the Kantorovich potential in~\eqref{eq:conjugate_pair}
satisfies~$\phi^c = - \phi$ and~$|\phi(x) - \phi(y)| \leq c(x,y)$ for all~$x,y \in \Omega$.
(b) Furthermore, the Kantorovich potential is differentiable almost everywhere in $\mathring{\Omega}$, with $|\nabla \phi | \leq \xi$~a.e.
\end{lemma}
%

It follows from Lemma~\ref{lemma:cconvex_phic} that
for absolutely continuous probability measures~$\mu, \mu^* \in \mathcal{P} \left( \Omega \right)$,
the optimal transport cost for the transport of~$\mu$ onto~$\mu^*$ can be
written as:
\begin{align}
  \begin{aligned}
    C(\mu, \mu^*) = \sup_{\phi \in W^{1,\infty}(\Omega, \mu)} &~\int_{\Omega}  \phi \left( 1 - \frac{d\mu^*}{d\mu} \right)~d\mu, \\
    							&\text{s.t.}~ |\nabla \phi | \leq \xi,~~\mu-\text{a.e.~in}~\Omega,
  \end{aligned}
  \label{eq:Kantorovich_duality_metric_cost}
\end{align}
where~$d\mu^* / d\mu$ is the Radon-Nikodym derivative.
With~$G_{\phi}(x) = \left[ | \nabla \phi(x) |^2 - \xi(x)^2 \right]/2$, 
we can rewrite the constraint in Problem~\eqref{eq:Kantorovich_duality_metric_cost} 
as $G_{\phi} \leq 0$, $\mu$--a.e. in~$\Omega$.
Furthermore, we have~$G_{\phi} \geq - \xi^2 / 2$ 
for any~$\phi \in W^{1,\infty}(\Omega, \mu)$.
Therefore, it follows that~$G_{\phi} \in L^{\infty}(\Omega, \mu)$, 
and the Lagrange multiplier corresponding to the constraint~$G_{\phi} \leq 0$ 
then belongs to~$L^{\infty}(\Omega, \mu)^*$, 
the dual space of~$L^{\infty}(\Omega, \mu)$.
Since $L^{\infty}(\Omega, \mu)^*$ is isomorphic~\cite{ND-JTS:58} to
$\mathrm{ba}(\Omega, \mu)$\footnote{The space of bounded finitely additive measures 
on~$\Omega$ that are absolutely continuous w.r.t.~$\mu$.},
we get that~$\lambda(G_{\phi}) = \frac{1}{2} \int_{\Omega} \left( |\nabla \phi|^2 - |\xi|^2 \right) \rho d\lambda$.
The Lagrangian functional $L: W^{1, \infty}(\Omega, \mu) \times
L^{\infty}(\Omega, \mu)^*_{\geq 0} \rightarrow \real$ can now be
defined as follows:
\begin{align*}
  L(\phi, \lambda) = - \int_{\Omega} \phi \left(1 -
    \frac{d\mu^*}{d\mu} \right) d\mu ~+~ \lambda(G_{\phi}).
\end{align*}
The following theorem establishes the existence of a (global) 
maximizer for Problem~\eqref{eq:Kantorovich_duality_metric_cost},
also called a Kantorovich potential,
and characterizes the saddle points of the Lagrangian~$L$:
\begin{theorem}[\bf \emph{First-order optimality conditions}]\label{thm:KKT_Kantorovich_duality}
Let $\mu, \mu^* \in \mathcal{P}(\Omega)$ be absolutely
continuous probability measures with densities $\rho, \rho^* \in L^1(\Omega)$.
Problem~\eqref{eq:Kantorovich_duality_metric_cost} has 
a global maximizer~$\phi_{\mu \rightarrow \mu^*} \in W^{1, \infty}(\Omega, \mu)$.
The Lagrangian~$L$ has a saddle point
$(\phi_{\mu \rightarrow \mu^*}, \lambda_{\mu \rightarrow \mu^*}) \in W^{1, \infty}(\Omega, \mu) \times L^1(\Omega, \mu)_{\geq 0}$. 
Moreover, $(\phi_{\mu \rightarrow \mu^*}, \lambda_{\mu \rightarrow \mu^*})$ satisfies
the first-order optimality conditions:
\begin{enumerate}
\item \textit{Stationarity:} 
The saddle point $(\phi_{\mu \rightarrow \mu^*}, \lambda_{\mu \rightarrow \mu^*})$ 
weakly satisfies the Poisson equation,
\begin{align*}
\begin{aligned}
	 - \frac{1}{\rho} \nabla \cdot \left( \rho \lambda_{\mu \rightarrow \mu^*} \nabla \phi_{\mu \rightarrow \mu^*}  \right) 
	 									&= 1 - \frac{\rho^*}{\rho},~~\mu-\text{a.e. in}~\Omega, \\
	 \rho \lambda_{\mu \rightarrow \mu^*} \nabla \phi_{\mu \rightarrow \mu^*} \cdot \mathbf{n} &= 0 \qquad \text{on}~\partial \Omega,
\end{aligned}
\end{align*}
where $\mathbf{n}$ is the outward normal to the boundary~$\partial \Omega$. \\
\item \textit{Feasibility:} 
$\lambda_{\mu \rightarrow \mu^*} \geq 0$ and $|\nabla \phi_{\mu \rightarrow \mu^*} | \leq  \xi$,~~$\mu$--a.e. in~$\Omega$. \\
\item \textit{Complementary slackness:} $\lambda_{\mu \rightarrow \mu^*} ( |\nabla \phi_{\mu \rightarrow \mu^*} | - \xi )= 0$,~$\mu$--a.e. in~$\Omega$.
\end{enumerate}
\end{theorem}

Having characterized the Kantorovich potential
via the saddle points of the Lagrangian~$L$,
in what follows we devise a stochastic process to achieve 
iterative optimal transport to the target measure~$\mu^*$.
We obtain a formulation of optimal transport of 
probability measures in which the net transport 
of a given initial probability measure~$\mu_0$ onto
a target probability measure~$\mu^*$ is
carried out by a stochastic process.
However, we first require the following lemma 
which establishes that the optimal 
transport can be decomposed into multiple stages:
\begin{lemma}[\bf \emph{Decomposition of optimal transport cost}]
\label{lemma:OT_stage_sequence}
Given atomless probability measures~$\mu, \mu^* \in \mathcal{P}(\Omega)$, the cost of optimal transport
from~$\mu$ to~$\mu^*$ satisfies:
\begin{align*}
	C(\mu, \mu^*) = \min_{ \nu \in \mathcal{P}(\Omega) }~C(\mu, \nu) + C(\nu, \mu^*).
\end{align*}
\end{lemma}

We note that Lemma~\ref{lemma:OT_stage_sequence} allows for the decomposition of optimal transport into an iterative process. By exploiting this decomposition, we will be able to integrate ``small transport steps'' with distributed multi-agent computations that approximate the Kantorovich potential of optimal transport, as explained in the following section.
Thus, as a first step, we evaluate decompositions for which 
the individual stage costs are upper bounded by
an~$\epsilon > 0$, as established by the  
following theorem:
\begin{theorem}[\bf \emph{Stochastic process for optimal iterative transport}]
\label{prop:law_Kant_duality_update_RV}
Let $\epsilon > 0$ and let
$\lbrace X(k) \rbrace_{k \in \mathbb{N}}$ be a stochastic
process with absolutely continuous marginals $\mu_k \in \mathcal{P}(\Omega)$ (i.e., $X(k) \sim \mu_k$ for every $k \in \mathbb{N}$)
such that for (measurable) $B \subseteq \Omega$ 
the following holds
\begin{align} \label{eq:Kantorovich_duality_update_RV}
\begin{aligned}
    \mathbb{P} \left[ X(k+1) \in  B \; | \; X(k) = x \right]
    = \frac{\ell \left( B \cap \mathcal{M}^{\epsilon}_{\mu_k}(x) \right)}{\ell \left(\mathcal{M}^{\epsilon}_{\mu_k}(x) \right)}
\end{aligned}
\end{align}
with the set $\mathcal{M}^{\epsilon}_\mu(x) = \arg \min_{z \in B^c_{\epsilon}(x)}~ c(x, z) + \phi_{\mu \rightarrow \mu^*}(z)$, where and $\ell$ denotes the arclength (or induced Lebesgue) measure.
Then, the sequence $\{ \mu_k \}_{k \in \mathbb{N}}$ satisfies 
\begin{align}
	\begin{aligned}
		\mu_{k+1} \in \arg \min_{\nu \in \mathcal{P}(\Omega)}~ &\left \lbrace C(\mu_k, \nu) + C(\nu, \mu^*), \right. \\
		&\quad \left. \text{s.t.}~ C(\mu_k, \nu) \leq \epsilon \right \rbrace,
	\end{aligned}
	\label{eq:OT_prob_iterative_scheme}
\end{align}
and $C(\mu_k, \mu^*) \rightarrow 0$, as $k \rightarrow \infty$. 

\end{theorem}

Some comments on Theorem~\ref{prop:law_Kant_duality_update_RV} are in order. 
We note that for a given $\mu \in \mathcal{P}(\Omega)$, at any $x \in \Omega$, the set of minimizers $\mathcal{M}_\mu^\epsilon(x)$ is the segment of the geodesic from $x$ to its image $T_{\mu \rightarrow \mu^*}(x)$ (the optimal transport map from $\mu$ to $\mu^*$)
that lies within the $\epsilon$ $c$-ball centered at $x$.
The update of the stochastic process involves uniformly sampling from this set. By limiting the update to distributions that are close in the sense of optimal transport, we achieve various objectives: first, we can account for the physical limitations of robots arising from the multi-agent setting wherein the agents cannot take arbitrarily large steps. Second, we can integrate such short displacement steps with distributed computations, as explained in the next section. Third, we also note that the resulting iterative process opens the door for feedback-based implementations, applicable to Model Predictive Control and
tracking time-varying~$\mu^*$. 
\section{Multi-agent implementation}
\label{sec:DOT}
\subsection{Multi-agent optimal transport algorithm}
We now obtain an implementable algorithm for multi-agent transport 
by the stochastic process~\eqref{eq:Kantorovich_duality_update_RV},
via a discretization of the Kantorovich potential.
At the culmination of the previous section, it becomes evident that while the process~\eqref{eq:Kantorovich_duality_update_RV} achieves
optimal iterative transport of absolutely continuous measures,
multi-agent systems are intrinsically represented by discrete measures,
whereby we are faced with the problem of obtaining an appropriate discretization and an implementable algorithm for the transport process.
To this end, we first discretize Problem~\eqref{eq:Kantorovich_duality_metric_cost}
onto a Voronoi partition generated by the positions of the agents
and devising a primal-dual algorithm to solve the discretized problem.
We then note that the primal-dual algorithm for the multi-agent 
system is Laplacian-based, thereby being intrinsically distributed
in nature, i.e., requiring communication only between immediate neighbors
on the Delaunay graph corresponding to the Voronoi partition.
This crucial fact enables scalable large-scale implementation of 
multi-agent optimal transport. In what follows, we first develop
the underlying structure for a discretization of 
Problem~\eqref{eq:Kantorovich_duality_metric_cost}, followed 
by the disctributed primal-dual algorithm.

Let~$\lbrace x_i(0) \rbrace_{i=1}^N$ be the positions of 
the~$N$ agents, distributed independently and identically 
according to a probability measure~$\mu_0$.
The idea is to transport the agents by the 
iterative process~\eqref{eq:Kantorovich_duality_update_RV}
to obtain~$\lbrace x_i(k) \rbrace_{i=1}^N$ at any time~$k$.
Let $\widehat{\mu}_N(k) = \frac{1}{N} \sum_{i=1}^N \delta_{x_i(k)}$ be the empirical 
measure generated by the agents~$\lbrace x_i(k) \rbrace_{i=1}^N$ at time~$k$.
To this end, we formulate a (finite) $N$-dimensional distributed optimization 
to be implemented by the agents to obtain local estimates of the Kantorovich potential.
We approximate the true Kantorovich potential by a~$\Phi^d : \naturals \times \Omega \rightarrow \real$
generated by an
(finite) $N$-dimensional vector~$\phi(k) = (\phi^1(k), \ldots, \phi^N(k)) \in \real^N$, 
such that~$\Phi^d(k, x_i(k)) = \phi^i(k)$ for~$i \in \lbrace 1, \ldots, N \rbrace$
and~$\Phi^d(k , x)$ for~$x \in \Omega \setminus \lbrace x_1(k), \ldots, x_N(k) \rbrace$ is defined by a suitable multivariate
interpolation.
In particular, let~$\lbrace \mathcal{V}_i(k) \rbrace_{i=1}^N$ be the Voronoi
partition of~$\Omega$ generated by~$\lbrace x_1(k), \ldots, x_N(k)
\rbrace$ w.r.t. the metric~$c$,
and~$\Phi^d = \sum_{i=1}^N \phi^{\mathcal{V}_i(k)}$
(decomposed into a sum of~$N$ functions~$\phi^{\mathcal{V}_i(k)}$ with
supports~$\mathcal{V}_i(k)$).  We assume that at time~$k$, the
agents~$i,j$ corresponding to neighboring cells~$\mathcal{V}_i(k)$
and~$\mathcal{V}_j(k)$ are connected by an edge, which defines a
connected graph \linebreak $G(k) = \left( \lbrace x_i(k)
  \rbrace_{i=1}^N, E(k) \right)$ (where~$E(k)$ is the edge set of the
graph~$G(k)$ at time~$k$).  

Dropping the index~$k$ (as is clear from
context), the finite dimensional approximation of the Kantorovich
duality~\eqref{eq:Kantorovich_duality_metric_cost} for the transport
between~$\widehat{\mu}_N$ and~$\mu^*$, restricted to the graph~$G$, is
given by:
\begin{align}
	\begin{aligned}
		\max_{(\phi^1, \ldots,\phi^N)}~ &\sum_{i=1}^N \left( \frac{1}{N} \cdot \phi^i  - \mathbb{E}_{\mu^*} [\phi^{\mathcal{V}_i}] \right) \\ \vspace*{0.1in}
		&\text{s.t.}~ |\phi^i - \phi^j| \leq c(x_i, x_j),~ \forall (i,j) \in E.
	\end{aligned}
	\label{eq:Kantorovich_dual_discrete}
\end{align}
We call~\eqref{eq:Kantorovich_dual_discrete} a restriction of~\eqref{eq:Kantorovich_duality_metric_cost} 
to the graph~$G$ because we only impose the constraint~$|\phi^i - \phi^j| \leq c(x_i, x_j)$ on neighbors~$i,j$ on the graph.

We note that Theorem~\ref{prop:law_Kant_duality_update_RV} guarantees convergence of the algorithm solving~\eqref{eq:Kantorovich_dual_discrete} exactly at each stage $k$ in the limit $N \rightarrow \infty$, when the closure of these samples are equal to the support of $\mu_k$. When $N$ is finite, there always will be an approximation error, which is expected to lead to an approximate convergence.  This is what it is verified in the simulations later.

We solve the optimization problem~\eqref{eq:Kantorovich_dual_discrete} by
a primal-dual algorithm, and its solution is used to update the agent positions
by~\eqref{eq:Kantorovich_duality_update_RV}. We take~$\Phi^d$ here
to be a simple function, such that~$\phi^{\mathcal{V}_i} (x) = \phi^i$
for~$x \in \mathcal{V}_i$. The Lagrangian for the problem~\eqref{eq:Kantorovich_dual_discrete},
with~$\Phi^d$ a simple function and~$c(x_i, x_j) = c_{ij}$, is given by:
\begin{align*}
\begin{aligned}
	L_d  = &- \sum_{i=1}^N  \phi^i \left[ \frac{1}{N}  - \mu^* ( \mathcal{V}_i) \right] \\
					&\qquad + \frac{1}{2} \sum_{i = 1}^N \sum_{j \in \mathcal{N}_i} \lambda_{ij} \left[ \left| \phi^i - \phi^j \right|^2 - c_{ij}^2 \right],
\end{aligned}
\end{align*}
and the primal-dual algorithm (with step size~$\tau$) is given~by (for $i \in \lbrace 1, \ldots, N \rbrace$
and $j \in \mathcal{N}_i$):
\begin{align}
	\begin{aligned}
		\phi^i(l+1) &= \phi^i(l) - \tau \sum_{j \in \mathcal{N}_i} \lambda_{ij}(l) \left( \phi^i(l) - \phi^j(l) \right) \\
							&\qquad \qquad \qquad + \left( \frac{1}{N} - \mu^*(\mathcal{V}_i) \right), \\
		\lambda_{ij} (l+1) &= \max \left \lbrace 0, \lambda_{ij}(l) + \tau \left[ \frac{1}{2} \left| \phi^i(l) - \phi^j(l) \right|^2 - c_{ij}^2 \right]  \right \rbrace.
	\end{aligned}
	\label{eq:p-d_Kantorovich_discrete}
\end{align}
The term $\sum_{j \in \mathcal{N}_i} \lambda_{ij}(l) \left( \phi^i(l)
    - \phi^j(l) \right)$ in~\eqref{eq:p-d_Kantorovich_discrete} is the
  action of the weighted Laplacian matrix (with
  weights~$\lambda_{ij}(l)$) on~$\phi(l)$.
  We note from the structure of the above algorithm that it renders itself to a
distributed implementation by the agents, where agent~$i$ uses
information from its neighbors~$j \in \mathcal{N}_i$ to
update~$\phi^i$ and~$\lbrace \lambda_{ij} \rbrace_{j \in
  \mathcal{N}_i}$. 

At the end of every step~$x_i (k) \mapsto x_i(k+1)$
from~\eqref{eq:Kantorovich_duality_update_RV}, the agent~$i$ 
assigns $\phi^i \leftarrow \Phi^d_{k} (x_i(k+1))$ as the initial
condition for the primal algorithm~\eqref{eq:p-d_Kantorovich_discrete} at
the time step~$k+1$ of the transport. 
Moreover, we are interested in an online
implementation of the transport, in that the agents do not wait for
convergence of the distributed primal-dual algorithm but carry out~$n$ iterations
of it for every update step~\eqref{eq:Kantorovich_duality_update_RV},
as outlined formally in the algorithm below.
\begin{algorithm}
\caption{Primal-dual based multi-agent optimal transport}
\textbf{Input:} Target measure $\mu^*$, Transport cost~$c(x,y)$, Bound on step size~$\epsilon$, Time step~$\tau$ \linebreak
\textbf{For each agent~$i$ at time instant~$k$ of transport:}
\begin{algorithmic}[1]
	\State Obtain: Positions~$x_j(k)$ of neighbors within communication/sensing radius~$r$ ($r \leq \text{diam}(\Omega)$, large enough to cover Voronoi neighbors)
	\State Compute: Voronoi cell~$\mathcal{V}_i(k)$, Mass of cell~$\mu^*(\mathcal{V}_i(k))$, Voronoi neighbors~$\mathcal{N}_i(k)$
	\State Initialize:~$\phi^i \leftarrow \Phi^d_{k-1} (x_i(k))$, $\lambda_{ij} \leftarrow \lambda_{ij}(k-1)$ (with~$\Phi^d_{0} = 0$, $\lambda_{ij}(0) = 0$)
	\State Implement~$n$ iterations of primal-dual algorithm~\eqref{eq:p-d_Kantorovich_discrete} (synchronously, in communication with neighbors~$j \in \mathcal{N}_i$)
				to obtain~$\phi^i(k)$,~$\lambda_{ij}(k)$
	\State Communicate with neighbors~$j \in \mathcal{N}_i$ to obtain~$\phi^j(k)$, construct local estimate of~$\Phi^d_{k}$ by multivariate interpolation
	\State Implement transport step~\eqref{eq:Kantorovich_duality_update_RV} with local estimate of~$\Phi^d_{k}$ (which approximates~$\phi_{\mu_k \rightarrow \mu^*}$)
\end{algorithmic}
	\label{alg:dist_opt_transport}
\end{algorithm}
\begin{remark}[\bf \emph{Optimize-then-discretize vs discretize-then-optimize}]
It can be seen that the key ingredient of the multi-agent optimal transport Algorithm~\ref{alg:dist_opt_transport} is the distributed computation by the agents of the (discretized) Kantorovich potential which is the maximizer in~\eqref{eq:Kantorovich_duality_metric_cost}.
This can be achieved either via a discretization of the 
PDE characterizing the first-order optimality condition
for the Kantorovich potential, namely the optimize-then-discretize
approach, or by first discretizing the infinite-dimensional optimization problem~\eqref{eq:Kantorovich_duality_metric_cost} 
over a partition (e.g., Voronoi) of the domain and then 
solving the finite-dimensional optimization problem resulting from the
discretization to directly obtain the discretized Kantorovich potential,
the latter being the discretize-then-optimize approach. 
In this paper we have
first laid the theoretical framework for the infinite-dimensional 
Kantorovich problem, which in principle represents the $N \rightarrow \infty$ limit. This allows for multi-agent algorithm design
by either approach. It is worth noting that while the 
distributed primal-dual algorithm~\eqref{eq:p-d_Kantorovich_discrete}
is presented as solving the discretized finite-dimensional optimization problem~\eqref{eq:Kantorovich_dual_discrete}, it 
can equivalently be viewed as computing the discretized solution of the 
PDE characterizing the first-order optimality in Theorem~\ref{thm:KKT_Kantorovich_duality}.
In the limits $N \rightarrow \infty$ of the number of agents and $n \rightarrow \infty$ of the number of iterations of the primal-dual algorithm, i.e. as the discretization more closely approximates the Kantorovich potential, we expect the two approaches to be equivalent, although this analysis is outside the scope of this paper.

\end{remark}
The primal-dual algorithm is scalable with respect to the number of agents, wherein the agents are only required to communicate with their spatial (Voronoi) neighbors, as seen from \eqref{eq:p-d_Kantorovich_discrete}. Voronoi partition-based distributed algorithms have been widely developed in the literature and scalable implementations exist~\cite{FB-JC-SM:09}. Complexity of internal iterations is explicitly controlled by the number of agents~$N$, the number of iterations~$n$ of the distributed primal-dual algorithm~\eqref{eq:p-d_Kantorovich_discrete} and computing $\mu^*(\mathcal{V}_i)$ across the Voronoi partition by the agents.
The runtime of the outer optimization is controlled by the distance between $\mu_0$ and $\mu^*$, the bound $\epsilon$ on the transport step size and the conformal factor $\xi$, while every agent incurs a cost to implement the optimization involved in the transport step that is distributed across the agents.
\subsection{Numerical experiments} \label{sec:simulation}
We now present numerical results for 
multi-agent optimal transport in~$\real^2$, based on the
the stochastic process~\eqref{eq:Kantorovich_duality_update_RV} (with~$c$
being the Euclidean metric and~$\epsilon = 0.02$),
where the local estimates of the Kantorovich potential are computed
by the distributed online algorithm~\eqref{eq:p-d_Kantorovich_discrete}
with a step size~$\tau = 1$, as outlined in Algorithm~\ref{alg:dist_opt_transport}. 

\begin{figure*}
\begin{center}
        \includegraphics[width=0.24\textwidth]{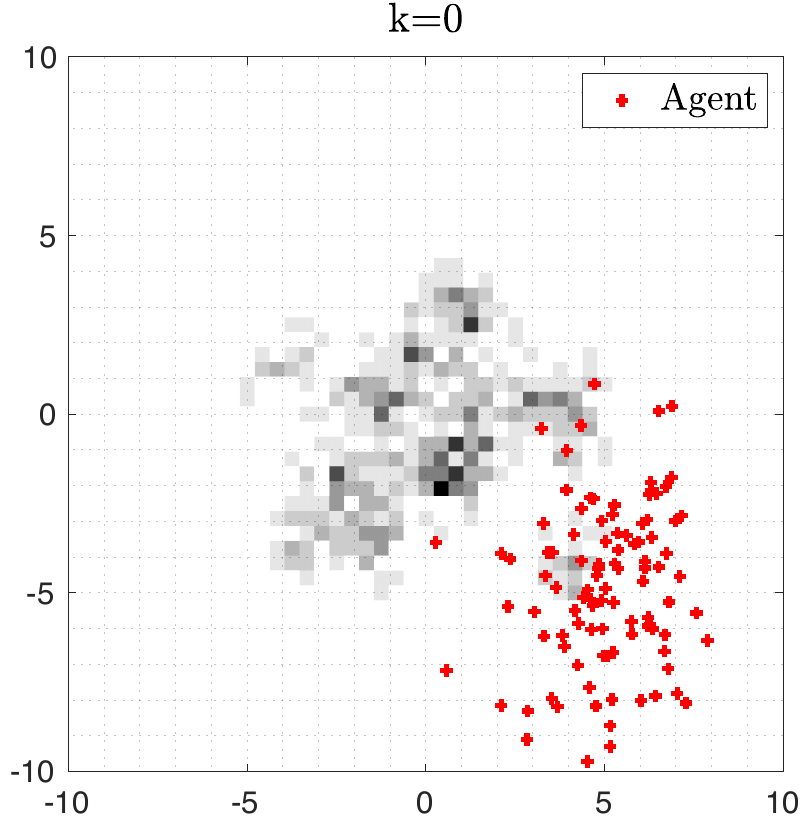}
                \includegraphics[width=0.24\textwidth]{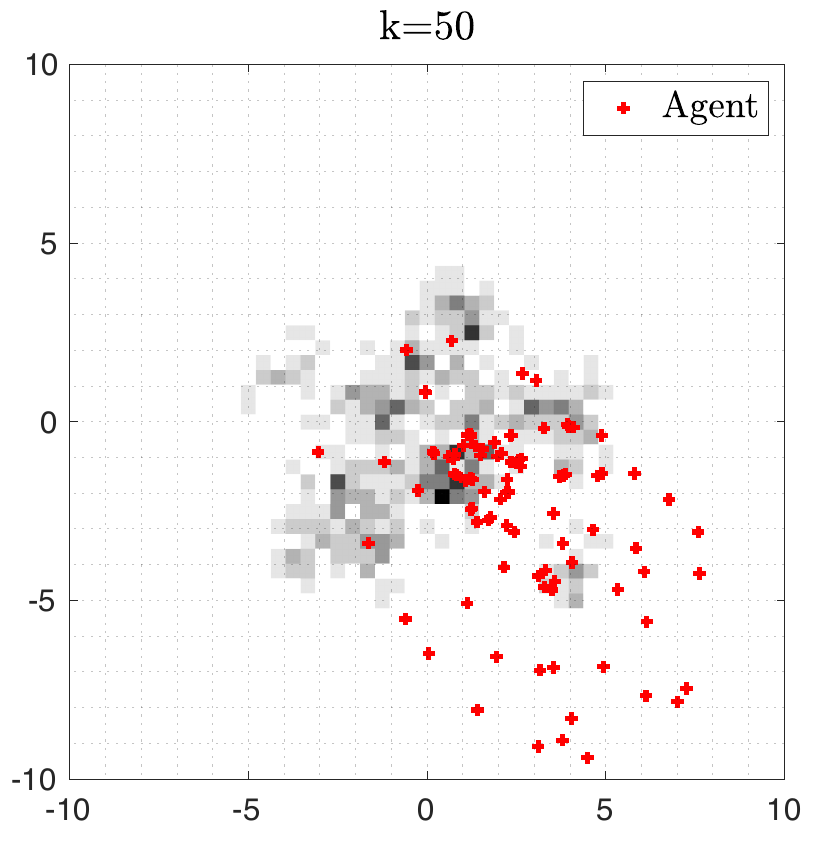}
                        \includegraphics[width=0.24\textwidth]{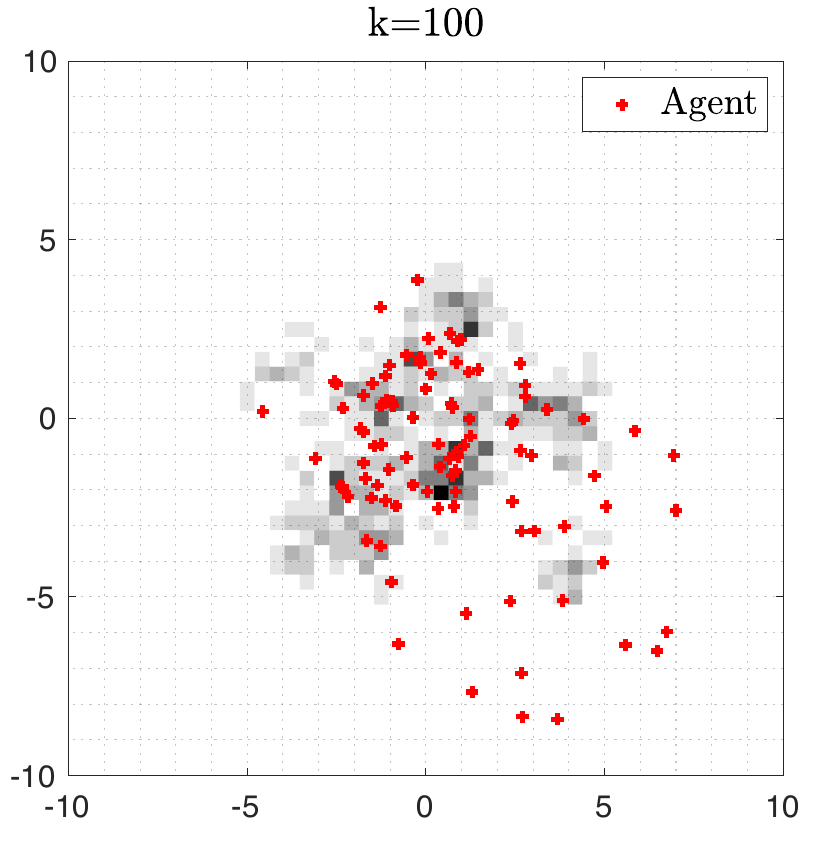}
                                \includegraphics[width=0.24\textwidth]{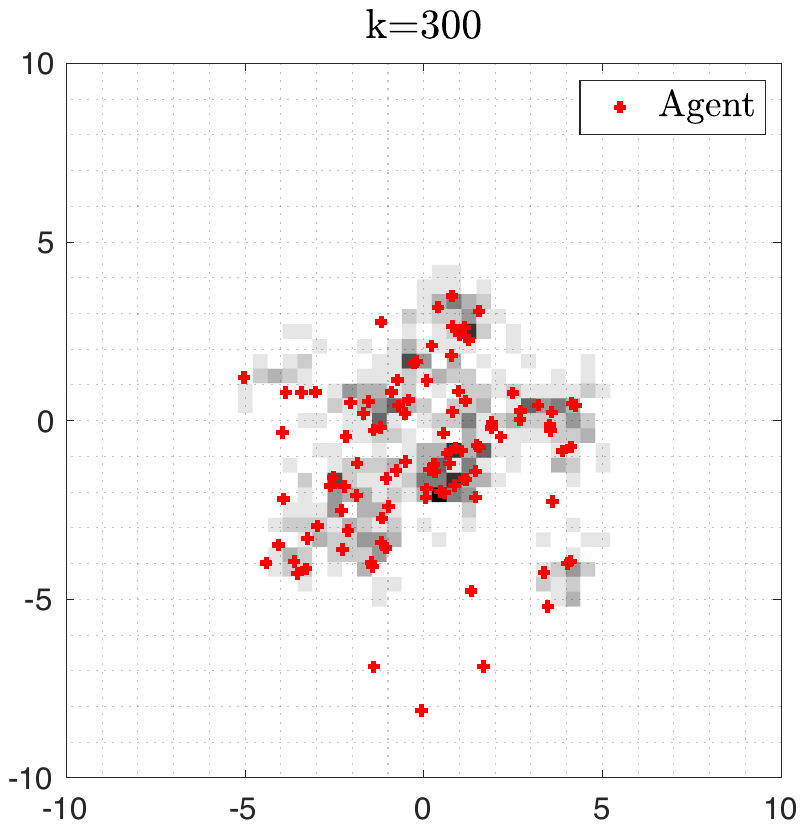}
\end{center}
        \caption{Positions of agents at
        				three different stages (time instants~$k=0, 50, 100, 300$) of transport by Algorithm~\ref{alg:dist_opt_transport},
       					i.e., the iterative process~\eqref{eq:Kantorovich_duality_update_RV} 
        				with local estimates of Kantorovich potential supplied by~\eqref{eq:p-d_Kantorovich_discrete} (with $n=1$ iterations of the primal-dual algorithm~\eqref{eq:p-d_Kantorovich_discrete} for every transport step~\eqref{eq:p-d_Kantorovich_discrete}; 
        				Target probability measure shown in grayscale with a darker shade indicating
          				a region of higher target density; The plots show convergence in time of the agents 
          				to full coverage of the target coverage profile (represented by the target probability distribution).}
	\label{fig:voronoi}
\end{figure*}

\begin{figure}[h]
    \includegraphics[width=0.45\textwidth]{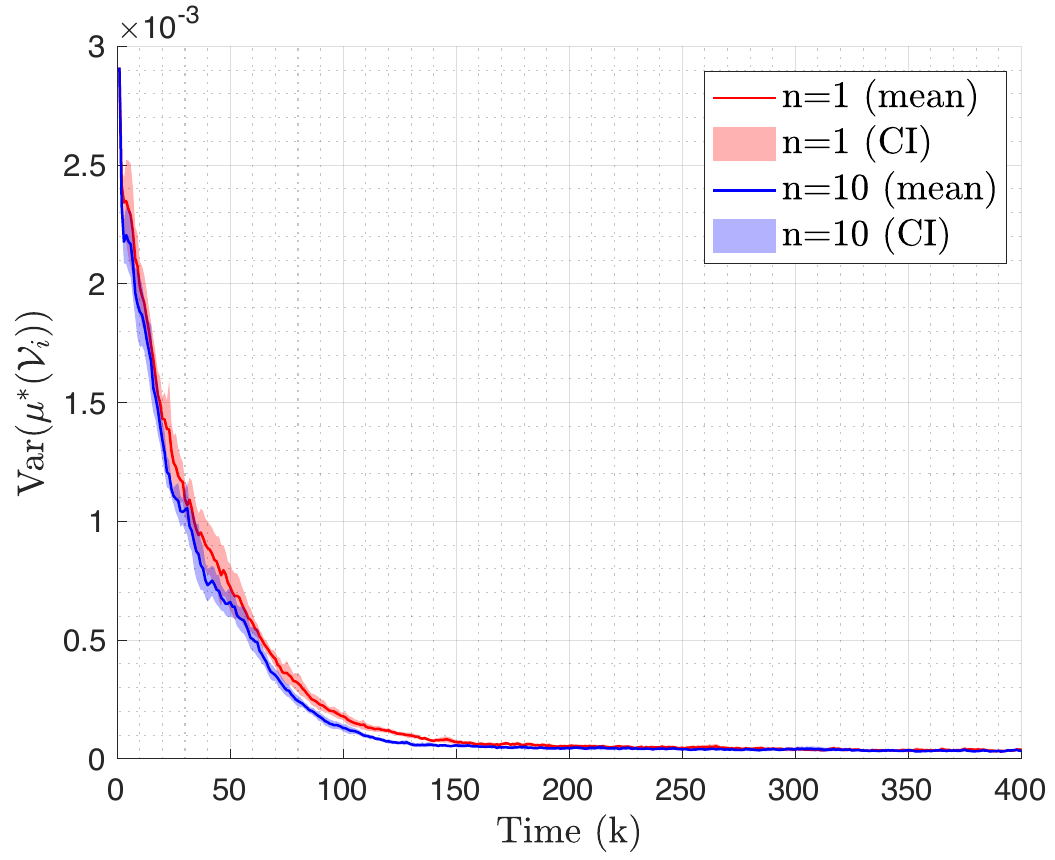}
	\caption{Variance in target mass~$\mu^*(\mathcal{V}_i)$ across the partition vs time
	for iteration steps~$n=1,10$ of the primal-dual algorithm~\eqref{eq:p-d_Kantorovich_discrete}
 for every transport step~\eqref{eq:Kantorovich_duality_update_RV}. The plot shows the
 (empirical) mean along with $95\%$ confidence bounds of the variance of~$\mu^*(\mathcal{V}_i)$
 from $10$ trials for each~$n$ from the same initial condition.
    }
	\label{fig:variance_plot}
\end{figure}

We chose as the target distribution the histogram of i.i.d. samples 
of a multimodal, bivariate Gaussian distribution (shown in grayscale in Figure~\ref{fig:voronoi}), 
and~$N = 100$ agents for the transport. Figure~\ref{fig:voronoi} shows the agents along with the
corresponding Voronoi partition of the domain, at three different
stages (time instants~$k=0, 50, 100, 300$) during the course of their
transport.  We observe that the agents are transported towards the
target probability measure and that a quantization of the target
measure is obtained. This is clarified further in
Figure~\ref{fig:variance_plot}, as described below.
As we had noted in the previous section, there exists a fundamental
trade-off between optimality and an online implementation of the
distributed optimal transport. We
sought to investigate the extent of this trade-off in simulation by
running multiple iterations~$n$ of the primal-dual algorithm~\eqref{eq:p-d_Kantorovich_discrete}
 for every iteration of the transport~\eqref{eq:Kantorovich_duality_update_RV}.  
The underlying rationale is that the distributed computation is many times faster
than the transport. Figure~\ref{fig:variance_plot} shows the rate of convergence
(w.r.t. the variance in target mass~$\mu^*(\mathcal{V}_i)$ across the partition) 
for~$n=1,10$. The randomness in Figure~\ref{fig:variance_plot} (as seen by the variation across trials and the $95\%$ confidence intervals) is partly due to the fact that the same initial distribution $\mu_0$ results in different samples across the trials, and the fact that the transport step~\eqref{eq:Kantorovich_duality_update_RV}  doesn't yield a unique minimizer but the agents are allowed to sample any one of the minimizers within an $\epsilon$-ball.
The numerical experiments were implemented in MATLAB R2023b on an Intel(R) Xeon(R) Platinum 8375C 2.90GHz CPU with $16$GB of RAM. The average runtime across $10$ independent runs of the experiment, for $k=400$ steps of the transport with $N=100$ agents, was $9.9$s. The average runtime for $n=20$ iterations of the distributed online primal-dual algorithm was $0.015$s.
%
%
%

\section{Conclusion}
\label{sec:conclusion}
We proposed a scalable, distributed iterative proximal
point algorithm for large-scale optimal transport of multi-agent
collectives. We obtained a dynamical formulation of optimal transport
of agents, for metric transport costs that are conformal to the
Euclidean distance.  We proposed a distributed primal-dual algorithm
to be implemented by the agents to obtain local estimates of the
Kantorovich potential, which are then used to implement the multi-agent
optimal transport. We studied the behavior of
the transport in simulation and investigated the
suboptimality of the online implementation. The analytic
characterization of the extent of the trade-off between optimality and the
online nature of the implementation is left for future~work.

\bibliographystyle{plain}
\bibliography{alias,SMD-add,SM,JC,FB}


\clearpage
\begin{appendices}

\section{Proofs of Lemmas}
\subsection{Proof of Lemma~\ref{lemma:cconvex_phic}} 
\label{app:proof_lemma_cconvex_phic}
{\emph{(a)}} This part of the proof employs similar arguments as the proof of the Kantorovich duality~\cite{CV:08}.
From~\eqref{eq:conjugate_pair}, we have:
\begin{align*}
  \phi(x) &= \inf_{y \in \Omega} \bigg( c(x,y) - \inf_{z \in \Omega} \left( c(z,y) - \phi(z) \right)  \bigg) \\
  &= \inf_{y \in \Omega} \sup_{z \in \Omega} \bigg( c(x,y) - c(z,y) + \phi(z)  \bigg) \\
  &\geq \inf_{y \in \Omega} \bigg( c(x,y) - c(z,y) + \phi(z)  \bigg) \\
  &= \inf_{y \in \Omega} \bigg( c(x,y) - c(z,y) \bigg) + \phi(z) \\ 
  &\geq - c(x,z) + \phi(z),
\end{align*}
where we have used the fact that~$c$ is a metric to obtain the final
inequality (for any~$y$, we have $c(x,y) - c(z,y) = c(x,y) - c(y,z) \geq - c(x,z)$, which
implies that $\inf_{y \in \Omega} \left(c(x,y) - c(z,y)\right) \geq -
c(x,z)$).  Moreover, since the above inequality holds for any~$x,z \in
\Omega$, we have~$|\phi(x) - \phi(z)| \leq c(x,z)$.

Now, when~$|\phi(x) - \phi(y)| \leq c(x,y)$, we have that $- \phi(x)
\leq c(x,y) - \phi(y)$, which implies that $ - \phi(x) \leq \inf_{y}
\left( c(x,y) - \phi(y) \right) = \phi^c(x)$.  Equivalently, we obtain
the relation $\phi(x) \geq - \phi^c(x)$.

  Similarly, from~\eqref{eq:conjugate_pair} $\phi^c(x) = \inf_{y} c(x,y) - \phi(y)$, we
  obtain $ \phi^c(x) \leq c(x,y) - \phi(y)$. By setting $y=x$ in the above inequality, 
  and using~$c(x,x) =0$ we get 
  $ \phi(x) \leq - \phi^c(x)$.
  In all, we have that~$\phi^c(x) = -\phi(x)$ and~$|\phi(x) -
  \phi(y)| \leq c(x,y)$.  

\emph{(b)}
We now have that~$\phi$ is Lipschitz continuous in~$\Omega$
(since~$c$ is conformal to the Euclidean metric from
and~$\Omega$ is compact). 
It follows from Rademacher's theorem\footnote{{\bf \emph{Rademacher's Theorem~\cite{GL:09}}:}
Let~$U \subset \real^d$ be open and bounded, 
and~$f : U \rightarrow \real^m$ be locally Lipschitz continuous in~$U$.
Then~$f$ is differentiable at almost every~$x \in U$.}
that~$\phi$ is differentiable almost everywhere in~$\Omega$. 
Also, from Theorems~4 and~5, Chapter~5.8 in~\cite{LCE:98}, we infer
that~$\phi \in W^{1, \infty}(\Omega)$ and that 
its (a.e.) gradient~$\nabla \phi$ is equal to its weak gradient
almost everywhere in~$\Omega$. 

For any minimizing curve $\gamma \in C^1([0,1]; \Omega)$ of \eqref{eq:conformal_metric}, we have 
\begin{align*}
    \lim_{h \rightarrow 0^+} \frac{c(\gamma(h),\gamma(0))}{h} 
    &= \lim_{h \rightarrow 0^+} \frac{1}{h} \int_0^h \xi(\gamma(t)) \; | \; \dot{\gamma}(t) | \dd t \\
    &= \xi(\gamma(0)) \left| \dot{\gamma}(0) \right|
\end{align*} 
Furthermore, note that for any $\psi \in C^1(\Omega)$ for which $|\psi(x) - \psi(y)|\le c(x,y)$ for all $x,y\in \Omega$, we get
\begin{align*}
   \left| \nabla \psi(\gamma(0)) \cdot \dot{\gamma}(0) \right| 
   &= \left| \left. \frac{d}{dh} \psi(\gamma(h)) \right|_{h=0} \right| \\
   &= \lim_{h \rightarrow 0^+} \frac{|\psi(\gamma(h)) - \psi(\gamma(0))|}{h} \\ 
   &\leq \lim_{h \rightarrow 0^+} \frac{c(\gamma(h),\gamma(0))}{h} \\
   &= \xi(\gamma(0)) \left| \dot{\gamma}(0) \right|.
\end{align*}
Since the above bound holds for any minimizing curve
$\gamma \in C^1([0,1]; \Omega)$ of \eqref{eq:conformal_metric},
it follows that 
$ \left| \nabla \psi(\gamma(0)) \right| \leq \xi(\gamma(0))$.
Now, applying the foregoing reasoning to 
the pointwise gradient of $\phi \in W^{1, \infty}(\Omega)$
satisfying $|\phi(x) - \phi(y)| \leq c(x,y)$ 
(which exists almost everywhere in $\Omega$ and is equal to its weak gradient),
we get that the weak gradient (which we interchangeably denote 
as $\nabla \phi$)
satisfies $ \left| \nabla \phi(x) \right| \leq \xi(x)$
almost everywhere in $\Omega$.  The above formal reasoning involved the use of a $C^1$ test function $\psi$ to define the derivative and carrying that over to the Lipschitz function $\phi$ wherever it is differentiable. 

\subsection{Proof of Lemma~\ref{lemma:OT_stage_sequence}} \label{app:proof_lemma_OT_stage_sequence}
The statement of Lemma~\ref{lemma:OT_stage_sequence} is a consequence of the metric property of~$C$ (Chapter~6,~\cite{CV:08}), which can be established using the triangle inequality for~$c$ along minimum-length geodesics within the Monge formulation~\eqref{eq:OT_Monge}.
From the Monge formulation~\eqref{eq:OT_Monge}, 
we have:
\begin{align*}
	C(\mu, \mu^*) = \inf_{ \substack{ T \in \mathcal{T}(\Omega) \\ T_{\#} \mu = \mu^* }} \int_{\Omega} c(x, T(x))~d\mu(x),
\end{align*}
where~$\mathcal{T}(\Omega)$ is the set of measureable maps over~$\Omega$.
For any measurable map~$T: \Omega \rightarrow \Omega$ such that $T_{\#} \mu_0 = \mu^*$
and any $x, z \in \Omega$, we have:
\begin{align*}
	c(x,T(x)) \leq c(x,z) + c(z,T(x)).
\end{align*}
Moreover, we have:
\begin{align*}
	c(x,T(x)) = \inf_{z \in \Omega} c(x,z) + c(z,T(x)),
\end{align*}
where the minimum is attained when~$z$ lies on the minimum-length geodesic 
from~$x$ to~$T(x)$ (the existence of minimum-length geodesics
follows from the Hopf-Rinow Theorem).
Since the above holds for any~$x \in \Omega$, it then follows that:
\begin{align*}
	&C(\mu, \mu^*) = \inf_{ \substack{ T \in \mathcal{T}(\Omega) \\ T_{\#} \mu = \mu^* }} \int_{\Omega} c(x,T(x))~d\mu(x)  \\ 
		&= \inf_{ \substack{ T \in \mathcal{T}(\Omega) \\ T_{\#} \mu = \mu^* }} \int_{\Omega} \inf_{z \in \Omega} \left[ c(x,z) + c(z,T(x)) \right] d\mu(x) \\
		&\stackrel{(a)}{=} \inf_{ \substack{ T \in \mathcal{T}(\Omega) \\ T_{\#} \mu = \mu^* }}  ~\inf_{ \substack{ \tilde{T} \in \mathcal{T}(\Omega) }}  \int_{\Omega} \left[ c(x, \tilde{T}(x)) + c(\tilde{T}(x),T(x)) \right] d\mu(x) \\
		&= \inf_{ \substack{ T \in \mathcal{T}(\Omega) \\ T_{\#} \mu = \mu^* }}  ~\inf_{ \substack{ \tilde{T} \in \mathcal{T}(\Omega) }}  \left[ \int_{\Omega} c(x, \tilde{T}(x)) d\mu(x) \right.  \\ & \qquad \qquad \qquad \qquad \left. +  \int_{\Omega} c(x,T \circ \tilde{T}^{-1} (x)) d \tilde{T}_{\#}\mu (x) \right] \\
		&=  \inf_{ \substack{ \tilde{T} \in \mathcal{T}(\Omega) }}  \left[ \int_{\Omega} c(x, \tilde{T}(x)) d\mu(x) \right. 
		\\ & \qquad \qquad \qquad \qquad \left. +  \inf_{ \substack{ T \in \mathcal{T}(\Omega) \\ T_{\#} \mu = \mu^* }} \int_{\Omega} c(x,T \circ \tilde{T}^{-1} (x)) d \tilde{T}_{\#}\mu (x) \right] \\
		&= \inf_{ \substack{\nu \in \mathcal{P}(\Omega) } } \left[ \inf_{ \substack{ \tilde{T}^{(1)} \in \mathcal{T}(\Omega) \\ \tilde{T}^{(1)}_{\#}\mu = \nu }} \int_{\Omega} c(x, \tilde{T}^{(1)}(x)) d\mu(x) \right. 
		\\ & \qquad \qquad \qquad \qquad \left. + \inf_{ \substack{ \tilde{T}^{(2)} \in \mathcal{T}(\Omega) \\ \tilde{T}^{(2)}_{\#}\nu = \mu^* }} \int_{\Omega} c(x, \tilde{T}^{(2)}(x)) d\nu(x)  \right] \\
		&= \inf_{ \substack{\nu \in \mathcal{P}(\Omega) } } C(\mu, \nu) + C(\nu, \mu^*),
\end{align*}
\noindent
where the equality~$(a)$ follows from an application of Theorem~3A in~\cite{RTR:76-b}.
We note that there clearly exists at least one minimizer~$\nu$ above 
(the choices~$\nu = \mu^*$ and $\nu = \mu$ trivially minimize the cost).
We therefore conclude that $C(\mu, \mu^*) = \min_{ \substack{\nu \in \mathcal{P}(\Omega) } } C(\mu, \nu) + C(\nu, \mu^*)$.

\section{Proofs of Theorems}
\subsection{Proof of Theorem~\ref{thm:KKT_Kantorovich_duality}} 
\label{app:proof_thm_KKT_Kantorovich_duality}
\emph{(i) Existence of Kantorovich potential~$\phi$:} We first rewrite
the Kantorovich duality~\eqref{eq:Kantorovich_duality_metric_cost} as
\begin{align*}
	C(\mu, \mu^*) = - \inf_{\phi \in W^{1, \infty}(\Omega, \mu)} ~&  \mathbb{E}_{\mu^*} [\phi] - \mathbb{E}_{\mu} [ \phi ], \\
							&\text{s.t.}~ \left| \nabla \phi \right| \leq \xi	,~\mu-\text{a.e.}
\end{align*} 
We note that the objective functional above is such that for any~$p
\in \real$ and any $\phi \in W^{1, \infty} (\Omega, \mu)$ with $\left|
  \nabla \phi \right| \leq \xi $, we have $- \mathbb{E}_{\mu} [ \phi +
p] + \mathbb{E}_{\mu^*} [\phi + p] = - \mathbb{E}_{\mu} [ \phi] +
\mathbb{E}_{\mu^*} [\phi]$.  Furthermore, the constraint is such that
$\left| \nabla (\phi + p) \right| = \left| \nabla \phi \right| \leq
\xi$. It follows that the objective functional and constraint in the
above optimization problem are invariant to constant vertical shifts
$\phi \mapsto \phi + p$.  In particular, for $\phi
\mapsto \phi - \mathbb{E}_{\mu}[\phi]$, we can rewrite the above as
\begin{align*}
	C(\mu, \mu^*) = - &\inf_{\phi \in W^{1, \infty}(\Omega, \mu)}  ~\underbrace{\mathbb{E}_{\mu^*} [\phi]}_{\triangleq L_{\mu^*}(\phi)}, \\
							&\text{s.t.}~ \begin{cases} &\mathbb{E}_{\mu} [ \phi ] = 0, \\
													&\left| \nabla \phi \right| \leq \xi, ~\mu-\text{a.e.} \end{cases}
\end{align*}
To see the invariance under the transformation, recall that the dual problem is given by
\begin{align*}
    C(\mu, \mu^*) = \sup_{\phi \in W^{1,\infty}(\Omega, \mu)} &~\int_{\Omega}  \phi \left( 1 - \frac{d\mu^*}{d\mu} \right)~d\mu, \\
    							&\text{s.t.}~ |\nabla \phi | \leq \xi,~~\mu-\text{a.e.~in}~\Omega,
\end{align*}
which, for absolutely continuous probability measures $\mu, \mu^* \in \mathcal{P}(\Omega)$, we can rewrite as
\begin{align*}
	C(\mu, \mu^*) = \sup_{\phi \in W^{1, \infty}(\Omega, \mu)} ~&  J(\phi) \triangleq \mathbb{E}_{\mu} [ \phi ] - \mathbb{E}_{\mu^*} [\phi], \\
							&\text{s.t.}~ \left| \nabla \phi \right| \leq \xi	,~\mu-\text{a.e.~in}~\Omega.
\end{align*} 
Let $\psi = \phi - \mathbb{E}_\mu [\phi]$. First see that
$\phi \in W^{1, \infty}(\Omega, \mu)$, if and only if $\psi \in W^{1, \infty}(\Omega, \mu)$. The objective function, evaluated at $\psi$, is given by
\begin{align*}
    J(\psi) &= \mathbb{E}_{\mu} [ \psi ] - \mathbb{E}_{\mu^*} [\psi] \\
&= \mathbb{E}_{\mu} [ \phi - \mathbb{E}_\mu [\phi] ] 
- \mathbb{E}_{\mu^*} [\phi - \mathbb{E}_\mu [\phi]] \\
&= \mathbb{E}_{\mu} [ \phi ] - \mathbb{E}_{\mu} [\mathbb{E}_\mu [\phi] ]
- \mathbb{E}_{\mu^*} [\phi] + \mathbb{E}_{\mu^*} [ \mathbb{E}_\mu [\phi]] \\
&= \mathbb{E}_{\mu} [ \phi ] - \mathbb{E}_{\mu^*} [\phi]
- \mathbb{E}_{\mu} [\mathbb{E}_\mu [\phi] ] + \mathbb{E}_{\mu^*} [ \mathbb{E}_\mu [\phi]]
\end{align*}
Observe that $\mathbb{E}_{\mu} [\mathbb{E}_\mu [\phi] ] 
= \mathbb{E}_{\mu^*} [ \mathbb{E}_\mu [\phi]] 
= \mathbb{E}_\mu [\phi]$ as $\mu$ and $\mu^*$ are probability measures. Therefore, it follows that
\begin{align*}
    J(\psi) &= \mathbb{E}_{\mu} [ \phi ] - \mathbb{E}_{\mu^*} [\phi] = J(\phi)
\end{align*}
Furthermore, note that $\nabla \psi = \nabla (\phi - \mathbb{E}_\mu [\phi])
= \nabla \phi - \nabla \left( \mathbb{E}_\mu [\phi] \right) = \nabla \phi$,
since $\nabla \left( \mathbb{E}_\mu [\phi] \right) = 0$ (this is due to the fact that 
$\mathbb{E}_\mu [\phi] = \int_\Omega \phi(x) d\mu(x)$, the average value of $\phi$ over $\Omega$, is a constant w.r.t $x \in \Omega$).

Note that $\mathcal{L}(\Omega, \mu) = \left \lbrace \phi \in W^{1,
    \infty}(\Omega, \mu) \; | \; \mathbb{E}_{\mu} [ \phi ] = 0,
\right.$ $\left. \left| \nabla \phi \right| \leq \xi, ~\mu-\text{a.e.}
\right \rbrace$ is closed, convex and bounded.  Boundedness
of~$\mathcal{L}(\Omega, \mu)$ follows from the compactness
of~$\mathbb{Y} = \cup_{\phi \in \mathcal{L}(\Omega, \mu)} \phi(\Omega)$ which implies that
there exists an~$M \in \realnonnegative$ such that $\mathbb{Y} \subset
B_M(0)$. It follows that for any~$\phi \in
\mathcal{L}(\Omega, \mu)$, we have~$\| |\phi| \|_{L^{\infty}(\Omega, \mu)} \leq
M$. Moreover, we have~$\| |\nabla \phi| \|_{L^{\infty}(\Omega, \mu)}
\leq \| |\xi| \|_{L^{\infty}(\Omega, \mu)}$.  Therefore, $\| \phi \|_{W^{1,\infty}(\Omega, \mu)} = \| |\phi| \|_{L^{\infty}(\Omega, \mu)} + \| |\nabla
\phi | \|_{L^{\infty}(\Omega, \mu)} \leq M + \| |\xi| \|_{L^{\infty}(\Omega, \mu)} < \infty$ for
any~$\phi \in \mathcal{L}(\Omega, \mu)$.

The objective functional $L_{\mu^*}$ is convex and lower semicontinuous 
(in fact, it is (Gateaux) differentiable~\cite{WR:64} as seen earlier
for absolutely continuous~$\mu^*$).
Let~$\lbrace \phi_n \rbrace_{n \in \mathbb{N}}$ be a minimizing sequence in~$\mathcal{L}(\Omega, \mu)$ 
for $L_{\mu^*}$, such that~$\phi_n \in \mathcal{L}(\Omega, \mu)$
and~$\lim_{n \rightarrow \infty} L_{\mu^*}(\phi_n)
= \inf_{\phi \in \mathcal{L}(\Omega, \mu)} L_{\mu^*}(\phi)$. Clearly,
the sequence~$\lbrace \phi_n \rbrace_{n \in \mathbb{N}}$ is uniformly
bounded since~$\| \phi_n \|_{W^{1,\infty}(\Omega, \mu)} \leq M + \sup_{x \in \Omega} \left| \xi \right|$.
It is also uniformly equicontinuous, since
$\left| \phi_n(x_1) - \phi_n(x_2) \right| \leq \left( \sup_{x \in \Omega} \left| \xi \right| \right) | x_1 - x_2 |$
for all~$n \in \mathbb{N}$. Therefore, by the Arzel{\`a}-Ascoli Theorem~\cite{WR:64}, 
there exists a uniformly convergent subsequence~$\lbrace \phi_{n_j} \rbrace_{j \in \mathbb{N}}$,
with the limit~$\phi^* \in \mathcal{L}(\Omega, \mu)$. 
Furthermore, by the continuity of~$L_{\mu^*}$,
we get~$\lim_{j \rightarrow \infty} L_{\sigma}(\phi_{n_j}) = L_{\mu^*}(\phi^*) 
= \min_{\phi \in \mathcal{L}(\Omega, \mu)} L_{\mu^*}(\phi)$.
By the convexity of the loss~$L_{\mu^*}$,
we get that~$\phi^*$ is a global minimizer of~$L_{\mu^*}$.

{\emph{(ii) Saddle points of Lagrangian functional~$L$:}}
The set $\mathcal{L}(\Omega)$ in the above problem can be expressed as
$\mathcal{L}(\Omega) = \left \lbrace \phi \in W^{1,\infty}(\Omega) \; | \; 
G(\phi) \in (-\infty , 0] \right \rbrace$, where the functional
$G(\phi) = \mathrm{ess}\sup |\nabla \phi|^2 - \xi^2$,
and we have the constraint qualification:
\begin{align*}
	0 \in  \mathrm{int} \left \lbrace G \left( W^{1,\infty}(\Omega) \right) - (-\infty , 0] \right \rbrace,
\end{align*}
where $G \left( W^{1,\infty}(\Omega) \right) - (-\infty , 0] = G \left( W^{1,\infty}(\Omega) \right) + [0,\infty)$,
and the operation $+$ denotes the Minkowski sum. 
This allows us to apply Theorem~3.6 in \cite{JFB-AS:13}
to infer that the set of Lagrange multipliers corresponding to the minimizer~$\bar{\phi}$
is a non-empty, convex, bounded and weakly$-^*$ compact subset of
$L^{\infty}(\Omega)^*_{\geq 0}$, and the set of Lagrange multipliers
is the same for any minimizer~$\bar{\phi}$.
Moreover, we note that $(-\infty, 0]$ is a closed convex cone, and it follows from
Theorem~3.4-(iii) in \cite{JFB-AS:13} that for any Lagrange multiplier~$\bar{\lambda}$,
the pair $(\bar{\phi}, \bar{\lambda})$ is a saddle point of the Lagrangian
functional~$L$.

We now evaluate the Gateaux derivative of $L_{\lambda}$\footnote{We adopt the notation $L_{\lambda}(\phi) = L(\phi, \lambda)$.} 
in the space $W^{1,\infty}(\Omega)$. For
  $\phi \in W^{1,\infty}(\Omega)$ and a variation $\eta \in
  W^{1,\infty}(\Omega)$, the directional
  derivative of $L_\lambda$ at $\phi$ along $\eta$ is given by:
  \begin{align*}
  \begin{aligned}
    &L_\lambda'(\phi,\eta) = \lim_{\epsilon \rightarrow 0^+}
    \frac{1}{\epsilon}\int_{\Omega} (\phi + \epsilon \eta -
      \phi)(-\rho + \rho^*)\dvol  \\  &\quad + \frac{1}{2\epsilon} \int_{\Omega} \left( |\nabla \phi +
      \epsilon \nabla \eta|^2 -|\xi|^2 - |\nabla \phi|^2 + |\xi|^2
      \right) d\lambda,
  \end{aligned}      
  \end{align*}
  which simplifies to:
   \begin{align*}
     L_\lambda'(\phi,\eta) &= \lim_{\epsilon \rightarrow 0^+}
     \frac{1}{\epsilon}\int_{\Omega} \epsilon \eta (-\rho +
       \rho^*) \dvol \\
       &\qquad \qquad + \frac{1}{2\epsilon} \int_{\Omega} \left( \epsilon^2 |\nabla \eta|^2  + 2
       \epsilon \nabla \phi \cdot \nabla \eta \right) d\lambda \\
     & = -\int_\Omega \eta (\rho - \rho^*) \dvol + \int_\Omega \nabla
     \phi \cdot \nabla \eta~ d\lambda ,
  \end{align*}
  where the dominated convergence theorem is applied in order to
  calculate the previous limit which exists for every $\eta \in
   W^{1,\infty}(\Omega)$.

%
%
 
  By the Minimax theorem and since~$(\bar{\phi},\bar{\lambda})$ 
  is a saddle point for the Lagrangian $L$, we have
  $L(\bar{\phi},\bar{\lambda}) = \sup_{\lambda} \inf_{\phi}L(\phi,\lambda) = \inf_{\phi}
  \sup_{\lambda} L(\phi,\lambda)$. Moreover, for any feasible~$\tilde{\phi}$ and $\lambda \geq 0$, 
  we note that $L(\tilde{\phi}, \lambda) \le L(\tilde{\phi} , 0) < +\infty$ and
  which implies that $\sup_{\lambda} L(\tilde{\phi},\lambda) \le L(\tilde{\phi},0)$.
  This implies $\sup_{\lambda} \inf_{\phi}L(\phi,\lambda) = 
  \inf_{\phi} \sup_{\lambda}L(\phi,\lambda) \le  L(\bar{\phi},0) < +\infty$. 
 Thus, for \emph{feasibility} 
  we must have $\int_{\Omega} (|\nabla \bar{\phi}|^2 - |\xi|^2)~d\bar{\lambda} \leq 0$,
  from which it follows that $|\nabla \bar{\phi}|^2 \le |\xi|^2$ a.e. in~$\Omega$
  as $\bar{\lambda} \in ba(\Omega)_{\geq 0}$ (since if
  $\bar{\phi}$ does not satisfy 
  $|\nabla \bar{\phi}|^2 \le |\xi|^2$ a.e. in~$\Omega$, 
  we would have that $\sup_\lambda L(\bar{\phi},\lambda) = + \infty$, 
  violating the finite upper bound). We also have:
  \begin{align*}
    &L(\bar{\phi},0) = - \int_\Omega \bar{\phi}(\rho - \rho^*)\dvol 
    \le L(\bar{\phi},\bar{\lambda}) 
    = \min_{\phi} L(\phi,\bar{\lambda}) \\
    &\le L(\bar{\phi}, \bar{\lambda})
    = -\int_\Omega \bar{\phi}(\rho - \rho^*)\dvol + \int_\Omega (|\nabla \bar{\phi}|^2 -  |\xi|^2)d\bar{\lambda} \\
    &\le -\int_\Omega \bar{\phi}(\rho - \rho^*)\dvol,
  \end{align*}
  where the first inequality is due to the saddle point definition, the
  second equality follows from the Minimax theorem, the last
  inequality from the fact that $\bar{\phi}$ is a feasible
  solution. All inequalities are indeed equalities, and we therefore
  see that \emph{complementary slackness} holds in the form
  $\int_\Omega (|\nabla \bar{\phi}|^2 - |\xi|^2)~d\bar{\lambda} = 0$.

Since the pair~$(\bar{\phi}, \bar{\lambda})$ is a saddle point of~$L$, it is also
a stationary point of~$L$ and we have~$L'_{\bar{\lambda}} (\bar{\phi}, \eta) = 0$ 
for any~$\eta \in W^{1,\infty}(\Omega)$,
and we get the \emph{stationarity} condition $\int_\Omega \eta (\rho - \rho^*) \dvol = \int_\Omega \nabla
 \bar{\phi} \cdot \nabla \eta~ d\bar{\lambda}$.

%
%
{\emph{(iii) Improved regularity of Lagrange multipliers:}}
We now establish stronger regularity for the Lagrange multipliers~$\bar{\lambda}$.
We have that the Lagrange multipliers $\bar{\lambda} \in L^{\infty}(\Omega)^*_{\geq 0}$,
which are finitely additive measures absolutely continuous w.r.t. the Lebesgue measure,
are also linear continuous functionals on~$L^{\infty}(\Omega)$ and must therefore
vanish on sets of Lebesgue measure zero (i.e., $\bar{\lambda}(A) = 0$ for
$A \subset \Omega$ with $\vol(A) = 0$).
Moreover, from Theorem~1.24 in~\cite{KY-EH:52},
we can decompose $\bar{\lambda} = \bar{\lambda}_c + \bar{\lambda}_p$,
where $\bar{\lambda}_c$ is a non-negative countably additive measure
and $\bar{\lambda}_p$ is non-negative and purely finitely additive.
By the Radon-Nikodym theorem, we get that there exists 
a function~$h_c \in L^1(\Omega)$ such that the countably additive
and absolutely continuous measure~$\lambda_c$ 
satisfies $d\bar{\lambda}_c = h_c \dvol$.
By substitution in the stationarity condition,
we get~$\int_\Omega \eta (\rho - \rho^*) \dvol = \int_\Omega (\nabla
     \bar{\phi} \cdot \nabla \eta)~h_c\dvol + \int_\Omega \nabla
     \bar{\phi} \cdot \nabla \eta~ d\bar{\lambda}_p$.
We now consider a set $D_{\delta} = \left \lbrace x \in \Omega \; | \;
 -\delta \leq |\nabla \phi (x)|^2 - \xi^2(x) \leq 0  \right \rbrace$,
 with $0 < \delta < \min_{x \in \Omega} \xi^2(x)$.
 By complementary slackness, we note that $\bar{\lambda}(\Omega 
 \setminus D_{\delta}) = 0$. Since $\bar{\lambda}_p$ is
 purely finitely additive, it implies that there must exist 
 a collection of nonempty sets $\lbrace E_n \rbrace_{N \in \mathbb{N}}$ with
 $E_{n+1} \subset E_n$ and $\lim_{n \rightarrow \infty} E_n = \emptyset$,
  such that $\lim_{n \rightarrow \infty} \bar{\lambda}_p(E_n) > 0$\footnote{
  For a countably additive measure~$\nu$ that is absolutely continuous w.r.t. the
  Lebesgue measure, and
  any collection of nonempty sets $\lbrace E_n \rbrace_{N \in \mathbb{N}}$ with
  $E_{n+1} \subset E_n$ and $\lim_{n \rightarrow \infty} E_n = \emptyset$,
  we have $\lim_{n \rightarrow \infty} \nu(E_n) = 0$~\cite{KY-EH:52}.}.
  Since $\bar{\lambda}(\Omega  \setminus D_{\delta}) = 0$, we can suppose
  without loss of generality that $E_0 \subset D_{\delta}$.
 We also consider another collection of nonempty sets $\lbrace E'_n \rbrace_{N \in \mathbb{N}}$,
 with the same properties (with $E'_0 \subset D_{\delta}$,
  $E'_{n+1} \subset E'_n$ and $\lim_{n \rightarrow \infty} E'_n = \emptyset$), 
 such that $E_n \subset E'_n$ for all $n \in \mathbb{N}$. 
 We note that for $x \in D_{\delta}$, we have $0 < \xi^2(x) - \delta \leq |\nabla \phi(x)|^2 \leq \xi^2(x)$,
 which implies that $\nabla \phi$ does not vanish on~$E'_n$ for any $n \in \mathbb{N}$.
 We now consider a family of variations $\eta_n \in W^{1,\infty}(\Omega)$ for $n \in \mathbb{N}$
 such that~$\eta_n$ and~$\nabla \eta_n$ are supported in~$E'_n$,
 $\nabla \phi \cdot \nabla \eta_n \geq 0$ in~$E'_n$
 and $\nabla \phi \cdot \nabla \eta_n \geq \epsilon$ in~$E_n$ (uniformly).
 The stationarity condition would now yield, for~$n \in \naturals$:
 \begin{align*}
 	&\int_{E'_n} \eta_n (\rho - \rho^*) \dvol \\
 	&= \int_{E'_n} (\nabla \bar{\phi} \cdot \nabla \eta_n)~h_c\dvol + \int_{E'_n} \nabla \bar{\phi} \cdot \nabla \eta_n~ d\bar{\lambda}_p \\
     &\geq  \int_{E'_n} (\nabla \bar{\phi} \cdot \nabla \eta_n)~h_c\dvol + \epsilon \int_{E_n} d\bar{\lambda}_p.
 \end{align*}
 In the limit $n \rightarrow 0$, we have $\lim_{n \rightarrow \infty} \int_{E'_n} \eta_n (\rho - \rho^*) \dvol = 0$
 and $\lim_{n \rightarrow \infty} \int_{E'_n} (\nabla
     \bar{\phi} \cdot \nabla \eta_n)~h_c\dvol = 0$, which 
     implies that $0 \leq  \lim_{n \rightarrow \infty} \epsilon \int_{E_n} d\bar{\lambda}_p
     \leq 0$, and we get $\lim_{n \rightarrow \infty} \bar{\lambda}_p(E_n) = 0$, i.e., 
     the measure~$\bar{\lambda}$ does not have a purely finitely additive component.
     Therefore, the measure~$\bar{\lambda}$ is countably additive (and absolutely continuous) 
     and possesses a Radon-Nikodym derivative w.r.t. the Lebesgue measure, in~$L^1(\Omega)$. 
     For ease of notation, we henceforth let $\bar{\lambda} \in L^1(\Omega)$ also denote its density function.
  Moreover, we note that since $\bar{\lambda} \in L^1(\Omega)_{\geq 0}$
  and $|\nabla \bar{\phi}|^2 - |\xi|^2 \leq 0$ a.e. in~$\Omega$,
  we can now indeed state the complementary slackness condition as 
   $\bar{\lambda}(|\nabla \bar{\phi} |^2 - |\xi|^2) = 0$ a.e. in~$\Omega$,
   which implies that $\bar{\lambda} (|\nabla \bar{\phi} | - |\xi|) =0$
  a.e. in~$\Omega$. We also see that~$\bar{\phi} \in W^{1, \infty} (\Omega)$.
 
  We recall from the stationarity condition that for any variation~$\eta \in W^{1,\infty}(\Omega)$ at~$\bar{\phi}$,
  we have $\int_\Omega \eta (\rho - \rho^*) \dvol = \int_\Omega \nabla
 \bar{\phi} \cdot \nabla \eta~ d\bar{\lambda} = \int_\Omega \bar{\lambda} \nabla
 \bar{\phi} \cdot \nabla \eta~\dvol$.
  Under stronger regularity of the saddle point~$(\bar{\phi}, \bar{\lambda})$, 
  the stationarity condition can be expressed as:
\begin{align*}
  & \int_{\Omega} (\rho - \rho^*) \eta~\dvol = \int_{\Omega} \bar{\lambda} \nabla \bar{\phi} \cdot \nabla \eta~\dvol \\
  &= - \int_{\Omega} \nabla \cdot (\bar{\lambda} \nabla \bar{\phi}) \eta~\dvol + \int_{\partial \Omega}
  \bar{\lambda} \nabla \bar{\phi} \cdot \mathbf{n} \eta~dS 
\end{align*}
where we have used the divergence theorem to obtain the final
equality, with~$S$ as the surface measure on~$\partial \Omega$. 
As the above holds for any variation~$\eta \in W^{1,\infty}(\Omega)$,
it must follow that~$-\nabla \cdot \left( \bar{\lambda} \nabla
  \bar{\phi} \right) = \rho - \rho^*$ in~$\Omega$ and~$\bar{\lambda}
\nabla \bar{\phi} \cdot \mathbf{n} = 0$ on~$\partial \Omega$,
and if we do not suppose stronger regularity of the saddle point~$(\bar{\phi}, \bar{\lambda})$,
the equations must be hold weakly.
Therefore, the saddle point~$(\bar{\phi}, \bar{\lambda})$
weakly satisfies~$-\nabla \cdot \left( \bar{\lambda} \nabla
  \bar{\phi} \right) = \rho - \rho^*$ in~$\Omega$ and~$\bar{\lambda}
\nabla \bar{\phi} \cdot \mathbf{n} = 0$ on~$\partial \Omega$.   

The above correspond to the necessary KKT conditions. Conversely, any
solution pair $(\bar{\phi},\bar{\lambda})$ which satisfies the above
KKT conditions results in a saddle point for the Lagrangian and it is
a solution to the original optimization problem.

\subsection{Proof of Theorem~\ref{prop:law_Kant_duality_update_RV}} \label{app:proof_prop_law_Kant_duality_update_RV}
{\emph{(a) Process~\eqref{eq:Kantorovich_duality_update_RV} achieves
\eqref{eq:OT_prob_iterative_scheme}}.}
We recall from~\eqref{eq:conjugate_pair} and
Lemma~\ref{lemma:cconvex_phic} that for the transport~$\mu_k \rightarrow \mu^*$,
and for any~$x \in \Omega$,
we have:
\begin{align}
	\phi_{\mu_k \rightarrow \mu^*}(x) = \inf_{y \in \Omega} c(x,y) + \phi_{\mu_k \rightarrow \mu^*}(y).
	\label{eq:Kant_pot_exp_min}
\end{align}
Clearly, for any~$x,y \in \Omega$, we have the inequality
$\phi_{\mu_k \rightarrow \mu^*}(x) \leq c(x,y) + \phi_{\mu_k \rightarrow \mu^*}(y)$.
Now, since~$c(x, y) = c(x,z) + c(z, y)$ for any (and only)~$z$ 
on a minimum-length geodesic from~$x$ to~$y$, we have
that the following holds $\mu_k$--a.e.:
\begin{align*}
	0 = \inf_{\substack{y \in \Omega, \\ z \in \bar{\gamma}_{x \rightarrow y}([0,1]) }} 
			~&\left[ - \phi_{\mu_k \rightarrow \mu^*}(x) + \phi_{\mu_k \rightarrow \mu^*}(z) - \phi_{\mu_k \rightarrow \mu^*}(z)
				\right. \\
				&\left. + \phi_{\mu_k \rightarrow \mu^*}(y)  + c(x, z) + c(z, y) \right],
\end{align*}
where~$\bar{\gamma}_{x \rightarrow y} : [0,1] \rightarrow \Omega$ is a length-minimizing geodesic from~$x$ to~$y$.
Rearranging the above, it follows that the following holds
$\mu_k$--a.e.:
\begin{align*}
	0 = \inf_{\substack{y \in \Omega, \\ z \in \bar{\gamma}_{x \rightarrow y}([0,1]) }} 
						~&\left[ - \phi_{\mu_k \rightarrow \mu^*}(x) + \phi_{\mu_k \rightarrow \mu^*}(z)  + c(x, z) \right] \\
				& + \left[ - \phi_{\mu_k \rightarrow \mu^*}(z) + \phi_{\mu_k \rightarrow \mu^*}(y)  + c(z, y) \right].
\end{align*}
Moreover, since the expressions (within the infimum) above are each non-negative, 
and their sum attains an infimum value of zero, 
we infer that they individually attain zero.
Thus, for any minimizer~$y^*$ in \eqref{eq:Kant_pot_exp_min}, we get that
every point on a length-minimizing geodesic from~$x$ to~$y^*$ is also a
minimizer. We note that the set of minimizers~$\mathcal{M}^\epsilon_{\mu_k}(x)$ in~\eqref{eq:Kantorovich_duality_update_RV} therefore
correspond to the segments of the length-minimizing geodesics (from the current iterate to the minimizers of~\eqref{eq:Kant_pot_exp_min}) 
contained in the $\epsilon$-ball $B^c_{\epsilon}$ centered at the current iterate.
Moreover, since $\phi_{\mu_k \rightarrow \mu^*}$ is differentiable $\mu_k$-almost everywhere,
it follows that the set of minimizers~$\mathcal{M}^\epsilon_{\mu_k}(x)$
are one-dimensional and are segments of geodesics from (almost every) $x \in \Omega$. 
It is known from Theorem~5.10 in~\cite{CV:08} that the optimal transport plans~$\pi_k \in \mathcal{P}(\Omega \times \Omega)$ 
(with~$\mu_k$ and~$\mu^*$ as the marginals) are supported in the set:
\begin{align*}
	\Gamma = \left \lbrace (x,y) \in \Omega \times \Omega \; | \; 
  \phi_{\mu_k \rightarrow \mu^*}(x) - \phi_{\mu_k \rightarrow \mu^*}(y)  = c(x,y) \right \rbrace.
\end{align*}
Let $\Gamma_x = \lbrace y \in \Omega \; | \; (x,y) \in \Gamma \rbrace$
and $\Gamma_x^{\epsilon} = \Gamma_x \cap B^c_{\epsilon}(x)$ for all $x \in \Omega$
where $\phi_{\mu_k \rightarrow \mu^*}$ is differentiable (since this is the case almost everywhere in $\Omega$,
we ignore the set of zero measure where $\phi_{\mu_k \rightarrow \mu^*}$ is not differentiable
in the rest of this proof). We see that
$\Gamma_{x}^{\epsilon}$ is the set of minimizers $\mathcal{M}^\epsilon_{\mu_k}(x)$ in~\eqref{eq:Kantorovich_duality_update_RV}.

Now, for any transport map~$\tilde{T}_k$ from~$\mu_k$ to~$\mu^*$,
we get $\phi_{\mu_k \rightarrow \mu^*}(x) \leq c(x,\tilde{T}_k(x)) + \phi_{\mu_k \rightarrow \mu^*}(\tilde{T}_k(x))$.
It then follows that:
\begin{align*}
	&\int_{\Omega} \left( \phi_{\mu_k \rightarrow \mu^*}(x) - \phi_{\mu_k \rightarrow \mu^*}(\tilde{T}_k(x)) \right) d\mu_k(x)  \\
	&= \int_{\Omega} \phi_{\mu_k \rightarrow \mu^*} d\mu_k - \int_{\Omega} \phi_{\mu_k \rightarrow \mu^*} d\mu^* \\
	&\leq \int_{\Omega} c(x, \tilde{T}_k(x)) d\mu_k(x).
\end{align*}
We see that the LHS is the optimal transport cost obtained from the Kantorovich dual formulation,
while an infimum over the RHS w.r.t.~$\tilde{T}_k$ would again yield the optimal transport cost from the
Monge formulation and an equality would then be attained. 
From Theorem~1 in~\cite{MF-RM:02}, we get that the Monge problem attains a
minimum, and let $\tilde{T}^*_k$ be an optimal transport map from~$\mu_k$ to~$\mu^*$.
Thus, we infer that $- \phi_{\mu_k \rightarrow \mu^*}(x) + \phi_{\mu_k \rightarrow \mu^*}(\tilde{T}^*_k(x)) + c(x, \tilde{T}^*_k(x)) = 0$,
$\mu_k$-almost everywhere in~$\Omega$, and $\tilde{T}^*_k(x) \in \Gamma_x$.
Moreover, the set $\Gamma_{x}^{\epsilon}$, which is
obtained as minimizers $\mathcal{M}^\epsilon_{\mu_k}(x)$ in~\eqref{eq:Kantorovich_duality_update_RV},
is precisely the segment of the geodesic (of length~$\epsilon$) from~$x$ to~$\tilde{T}^*_k(x)$
for $\mu_k$-almost every~$x \in \Omega$. From this correspondence between the Monge problem and the process~\eqref{eq:Kantorovich_duality_update_RV},
and from the the proof of Lemma~\ref{lemma:OT_stage_sequence}, 
it follows that the optimal transport cost $C(\mu_k, \mu^*)$ can be
decomposed as in~\eqref{eq:OT_prob_iterative_scheme} by the process~\eqref{eq:Kantorovich_duality_update_RV}.

{\emph{(b) Process \eqref{eq:Kantorovich_duality_update_RV} achieves
convergence $C(\mu_k, \mu^*) \rightarrow 0$}.}
We now consider a single update step of the process $\lbrace X(k) \rbrace_{k \in \mathbb{N}}$. 
For simplicity of notation, let $X$ be a random variable such that 
$X \sim \mu$ for some $\mu \in \mathcal{P}(\Omega)$ and let $x \in \Omega$ be a sample of $X$. 
Let $x^+$ be a sample drawn uniformly at random from the set of minimizers
$\arg \min_{z \in B^c_{\epsilon}(x)}~ c(x, z) + \phi_{\mu \rightarrow \mu^*}(z) = \arg \min_{z \in B^c_{\epsilon}(x)}~ c(x, z) + c(z, T_{\mu \rightarrow \mu^*}(x) )$.
Note that any such~$x^+$ satisfies $c(x^+, T_{\mu \rightarrow \mu^*}(x)) \leq c(x, T_{\mu \rightarrow \mu^*}(x))$.
Since this holds for any sample $x$ of $X$, we get that 
the random variable $X^+$ obtained from the single update 
of $X$, with $X^+ \sim \mu^+$, is such that 
 $C(\mu^+, \mu^*) \leq C(\mu, \mu^*)$. 
Therefore, the sequence of measures $\lbrace \mu_k \rbrace_{k \in \mathbb{N}}$ corresponding to the process $\{ X_k \}_{k \in \mathbb{N}}$ satisfies $C(\mu_{k+1}, \mu^*) \leq C(\mu_k, \mu^*)$
 for all $k \in \mathbb{N}$, and applying the monotone convergence
 theorem, we get that $C(\mu_k, \mu^*) \rightarrow \bar{C}$.
 By Prokhorov's Theorem~\cite{CV:08}, from which we get that $\mathcal{P}(\Omega)$ is compact w.r.t. the topology of weak convergence, and the fact that $\mu_k$ lies on the (unique)
geodesic between $\mu_0$ and $\mu^*$, it follows that 
there exists an accumulation point $\bar{\mu}$ on the
geodesic between $\mu_0$ and $\mu^*$ such that $\mu_k \rightarrow \bar{\mu}$.
This implies that $C(\bar{\mu}, \bar{\mu}^+) = 0$, where $\bar{\mu}^+ \in \mathcal{P}(\Omega)$ is obtained from a one step update of $\bar{X} \sim \bar{\mu}$. Furthermore, we have
\begin{align*}
    C(\bar{\mu}, \bar{\mu}^+) = \min_{\sigma \in \Pi(\bar{\mu}, \bar{\mu}^+)} \int_{\Omega \times \Omega} c(x,z) d\sigma(x,z).
\end{align*}
From Section~(a) of the proof, it follows that the process~\eqref{eq:Kantorovich_duality_update_RV}
results in optimal transport from $\bar{\mu}$ to $\bar{\mu}^+$.
With $\pi$ as the transition probability for the process~\eqref{eq:Kantorovich_duality_update_RV}, we get
\begin{align*}
    C(\bar{\mu}, \bar{\mu}^+) 
    &= \min_{\sigma \in \Pi(\bar{\mu}, \bar{\mu}^+)} \int_{\Omega \times \Omega} c(x,z) d\sigma(x,z) \\
    &= \int_{x \in \Omega} \int_{z \in \Omega} c(x,z) d\pi(z \; | \; x) d\bar{\mu}(x) \\
    &= \int_{x \in \Omega} \frac{1}{\ell(\mathcal{M}^{\epsilon}_{\bar{\mu}}(x))} \int_{z \in \mathcal{M}^{\epsilon}_{\bar{\mu}}(x)} c(x,z) d\ell(z) d\bar{\mu}(x).
\end{align*}
Suppose that $C(\bar{\mu}, \mu^*) = \eta > 0$. 
First note that for every $x \in \Omega$, the following holds
\begin{align*}
    c(x, T_{\mu \rightarrow \mu^*}(x)) = \min_{z \in B^c_{\epsilon}(x)}~ c(x, z) + c(z, T_{\mu \rightarrow \mu^*}(x) ),
\end{align*}
owing to the triangle inequality being an equality along the $c$-geodesic.
Furthermore, note that for any minimizer $z$ above, i.e. for $z \in \mathcal{M}^{\epsilon}_{\bar{\mu}}(x)$ 
we have $c(z, T_{\mu \rightarrow \mu^*}(x) ) \leq c(x, T_{\mu \rightarrow \mu^*}(x))$.
We then have
\begin{align*}
    \eta &= C(\bar{\mu}, \mu^*) = \int_{\Omega} c(x,T_{\bar{\mu} \rightarrow \mu^*}(x)) ~d\bar{\mu}(x) \\
    &= \int_{\Omega} \min_{z \in B^c_{\epsilon}(x)} \left[c(x, z) + c(z, T_{\mu \rightarrow \mu^*}(x) ) \right] ~d\bar{\mu}(x) \\
    &= \int_{\Omega} \frac{1}{\ell(\mathcal{M}^{\epsilon}_{\bar{\mu}}(x))} \int_{z \in \mathcal{M}^{\epsilon}{\bar{\mu}}(x)}
    \left[c(x, z)  \right. \\
    &\qquad \qquad \left. + c(z, T_{\mu \rightarrow \mu^*}(x) ) \right]  d\ell(z) d\bar{\mu}(x) \\
    &= C(\bar{\mu}, \bar{\mu}^+) \\
    &+ \int_{\Omega} \frac{1}{\ell(\mathcal{M}^{\epsilon}_{\bar{\mu}}(x))} \int_{z \in \mathcal{M}^{\epsilon}_{\bar{\mu}}(x)}
    c(z, T_{\mu \rightarrow \mu^*}(x) ) d\ell(z) d\bar{\mu}(x)
\end{align*}
From $C(\bar{\mu}, \bar{\mu}^+) = 0$ (since $\bar{\mu}$ is 
an accumulation point), we get 
\begin{align*}
    \int_{\Omega} \frac{1}{\ell(\mathcal{M}^{\epsilon}_{\bar{\mu}}(x))} \int_{z \in \mathcal{M}^{\epsilon}_{\bar{\mu}}(x)}
    c(z, T_{\mu \rightarrow \mu^*}(x) ) d\ell(z) d\bar{\mu}(x) = \eta 
\end{align*}
which implies that 
$\bar{\mu}$-a.e. in $\Omega$
\begin{align*}
    c(x,T_{\bar{\mu} \rightarrow \mu^*}(x))
    = \frac{1}{\ell(\mathcal{M}^{\epsilon}_{\bar{\mu}}(x))} \int_{z \in \mathcal{M}^{\epsilon}_{\bar{\mu}}(x)}
    c(z, T_{\mu \rightarrow \mu^*}(x) ) d\ell(z).
\end{align*}
Since for $z \in \mathcal{M}^{\epsilon}_{\bar{\mu}}(x)$ 
we have $c(z, T_{\mu \rightarrow \mu^*}(x) ) \leq c(x, T_{\mu \rightarrow \mu^*}(x))$, it follows that $\mathcal{M}^{\epsilon}_{\bar{\mu}}(x) = \{x\}$,
$\bar{\mu}$-a.e. in $\Omega$, which is the case
if and only if $T_{\bar{\mu} \rightarrow \mu^*}(x) = x$,
i.e., $C(\bar{\mu}, \mu^*) = 0$, and we get a contradiction.
Therefore, it follows that $\mu^*$ is the only accumulation point
and the sequence $\{ \mu_k \}_{k \in \mathbb{N}}$ converges to $\mu^*$
as $k \rightarrow \infty$ w.r.t. the topology of weak convergence in $\mathcal{P}(\Omega)$, i.e., $C(\mu_k, \mu^*) \rightarrow 0$.


\end{appendices}

\end{document}